\documentclass[11pt,a4paper]{amsart}
\usepackage{tikz}
\usepackage{amsmath,amsfonts,xargs,microtype,amssymb,bbm}
\microtypesetup{expansion=false}
\usepackage[utf8]{inputenc}
\usepackage{graphicx} 
\usepackage{amsthm}
\usepackage{pgfplots}
\pgfplotsset{compat=1.18}

\definecolor{mygreen}{RGB}{56,132,72}
\usepackage{tikz}
\usepackage{enumitem}
\usetikzlibrary{calc}
\usetikzlibrary{arrows.meta}
\usepackage{xfp}
\usepackage{tkz-tab}
\usepackage{theoremref,float}
\usepackage{thmtools, thm-restate}
\declaretheorem{theorem}
\usepackage{thm-restate}
\usepackage{hyperref}
\usepackage{comment}
\usepackage{color}
\usepackage{mathtools}
\usepackage{subcaption}
\usepackage{hyphenat}
\mathtoolsset{showonlyrefs}

\newtheorem{fact}{Fact}
\newtheorem{claim}{Claim}

\hypersetup{colorlinks=true,urlcolor=black, pdftitle="Convexification"}

\begin{document}

\newcommand{\R}{ {\mathbb R}}
\newcommand{\N}{ {\mathbb N}}
\newcommand{\s}{ {\mathbb S}}
\newcommand{\conv}{{\operatorname{conv}}}
\newcommand{\Span}{{\operatorname{span}}}
\newcommand{\diag}{{\operatorname{diag}}}
\newcommand{\Int}{{\operatorname{int}}}
\newcommand\numberthis{\addtocounter{equation}{1}\tag{\theequation}}
\newcommand{\footnotesmall}[1]{\footnote{{\fontsize{6.5pt}{8pt}\selectfont #1}}}

\theoremstyle{plain}

\newcommand{\LqPicture}[2]{%
\begin{scope}[yshift=#2]

\pgfmathsetmacro{\p}{#1}
\pgfmathsetmacro{\q}{\p/(\p-1)}
\pgfmathsetmacro{\seg}{pow(2,1/\p)-1}

\begin{axis}[
    width=7cm,
    height=7cm,
    xmin=-0.08, xmax=2.1,
    ymin=-0.08, ymax=2.1,
    axis lines=middle,
    axis line style={black, thick},
    xtick={0,1},
    ytick={0,1},
    xticklabels={$0$,$1$},
    yticklabels={$0$,$1$},
    tick style={draw=none},
    clip=false,
    enlargelimits=false,
    axis equal image,
    xlabel={},
    ylabel={},
    scale = .5
]

\addplot[
    domain=0:90,
    samples=300,
    very thick,
    mygreen
]
({pow(abs(cos(x)),2/\q)}, {pow(abs(sin(x)),2/\q)});

\addplot[
    very thick,
    mygreen
]
coordinates {(1,0) ({1+\seg},0)};

\addplot[
    very thick,
    mygreen
]
coordinates {(0,1) (0,{1+\seg})};

\addplot[
    only marks,
    mark=*,
    mark size=2pt,
    mygreen
]
coordinates {
    (0,1)
    (1,0)
    ({1+\seg},0)
    (0,{1+\seg})
};

\pgfmathtruncatemacro{\qden}{#1-1}
\node at (axis cs:0.68,-0.5)
    {$p=#1,\ q=\fpeval{round(#1/(#1-1),3)}$};
\end{axis}

\end{scope}
}

\newtheorem{lemma}[theorem]{Lemma}
\newtheorem{corollary}[theorem]{Corollary}
\newtheorem{prop}[theorem]{Proposition}

\newtheorem{conj}{Conjecture}

\theoremstyle{remark}
\newtheorem{remark}{Remark}

\makeatother

\title{Convexification effects of $L_p$-sums}

%

\author[M. Fradelizi]{Matthieu Fradelizi}
\address[M. Fradelizi]{LAMA, Univ Gustave Eiffel, Univ Paris Est Creteil, 77447 Marne-la-Vall\'ee, France.}
\email{matthieu.fradelizi@univ-eiffel.fr}

\author[A. Manui]{Auttawich Manui}
\address[A. Manui]{Department of Mathematics, Syracuse University, Syracuse, NY 13244, USA}
\email{amanui@syr.edu}

\author[M. Meyer]{Mark Meyer}
\address[M. Meyer]{Institut de Math\'ematiques de Jussieu, Sorbonne Universit\'e, CNRS, 4 Place Jussieu, 75252 Paris, France}
\email{meyer@imj-prg.fr}

\author[M. Szusterman]{Maud Szusterman}
\address[M. Szusterman]{Centre de Math\'ematiques Laurent Schwartz (CMLS), \'Ecole Polytechnique, 91128 Palaiseau Cedex, France}
\email{maud.szusterman@polytechnique.edu}

\author[A. Zvavitch]{Artem Zvavitch}
\address[A. Zvavitch]{Department of Mathematical Sciences, Kent State University, Kent, OH 44242, USA}
\email{zvavitch@math.kent.edu}


\date{}

\begin{abstract}
It is a classical result, quantified by Shapley, Folkman, and Starr, that the Minkowski averages of any compact set asymptotically approach its convex hull in the Hausdorff metric.
In this paper, we investigate the analogous convexification phenomenon for the Firey $L_p$-sums of compact sets. First, we establish new structural properties of the $L_p$-sum. Notably, we show that the iterated $L_p$-sum of compact sets whose convex hull contains the origin is eventually starlike. We then study $L_p$-convexification for three classical measures of nonconvexity: Hausdorff metric, volume deficit, and effective standard deviation. For each of these measures, we prove that the $k$-fold $L_p$-sum averages of a compact set $A$ converge to $\textup{conv}(0,A)$ when $1<p<\infty$, and to $\textup{conv}(A)$ when $p=\infty$. In each case, we show that the convergence to the convex hull is eventually monotone. For the Hausdorff metric and effective standard deviation, we prove sharp upper bounds for the $L_p$-sum of two compact sets.
\end{abstract}

\maketitle

\section{Introduction}
The Minkowski sum of nonempty sets $A,B\subset \mathbb{R}^n$ is defined by
\begin{align*}
    A+B:=\{a+b:a\in A,b\in B\}.
\end{align*}
Minkowski summation has a well-known convexification effect. Even though a set $A\subset\mathbb{R}^n$ may not be convex, it has been observed that its averages $k^{-1}\sum_{i=1}^kA$ converge to the convex hull of $A$ when $k$ goes to infinity. This phenomenon was quantified in terms of the Hausdorff metric by Shapley, Folkman, and Starr \cite{ST-69}, and by Emerson and Greenleaf \cite{EG-1969}. There has since been an increasing body of research studying the convexification effect from different perspectives. For more information on this topic, we direct the reader to the survey \cite{FMMZ-18}.

The objective of this paper is to study the convexification effect of the Firey $L_p$-sum \cite{F-62}, which is defined as follows. Given a convex body $K\subset\R^n$ (i.e., a convex, compact set with nonempty interior), the function
\begin{align*}
    h_K(u):=\sup_{x\in K}\langle x,u\rangle,\qquad u\in S^{n-1},
\end{align*}
is called the support function of $K$. If $K,L\subset \R^n$ are convex bodies that contain the origin, and if $1\leq p<\infty$, then the Firey $L_p$-sum of $K$ and $L$ is the set $K\oplus_p L$ whose support function satisfies 
\begin{align}\label{eq:lp_sum_support_function_def}
    h_{K\oplus_p L}(u)=(h_K(u)^p+h_L(u)^p)^{1/p}.
\end{align}
Taking the limit as $p\rightarrow \infty$, one arrives at the natural definition for the $L_{\infty}$-sum 
\begin{align*}
    h_{K\oplus_{\infty} L}(u):=\max\{h_K(u),h_L(u)\},
\end{align*} 
so that $K\oplus_{\infty} L=\text{conv}(K,L)$.
In the paper \cite{LYZ-12}, Lutwak, Yang, and Zhang extended the definition of Firey summation to the nonconvex setting, showing that if $A,B\subset\R^n$ are compact not necessarily containing the origin, then the definition 
\begin{align}\label{eq:lp_sum_union_def}
    A\oplus_pB:=\{t^{1/q}a+(1-t)^{1/q}b:0\leq t\leq 1,a\in A,b\in B\},\qquad \frac{1}{p}+\frac{1}{q}=1
\end{align}
is a natural extension of definition \eqref{eq:lp_sum_support_function_def}. In the special case $p=1$, definition \eqref{eq:lp_sum_union_def} is the usual Minkowski sum. Notice that in general, $A\oplus_\infty B\subsetneq\conv(A,B)$.

Motivated by \cite{FMMZ-18}, we will show that Firey $L_p$-sums have both starshaping and convexification effects. We detail our contributions below, which demonstrate these effects from the perspective of three measures: volume, Hausdorff distance from the convex hull, and effective standard deviation.

\subsection{Starshaping effect of \texorpdfstring{$L_p$}{Lp}-sums.}
In this paper, we call a set starshaped when it is starshaped at the origin. We write $A^{\oplus_pk}:=A\oplus_p\dots\oplus_pA$ for the $k$-fold $L_p$-sum of $A$ with itself. We prove that if $0\in \textup{conv}(A)$, then $A^{\oplus_p k}$ is starshaped for $k\ge n+1$. More generally, using the colorful Carath\'eodory theorems by B\'ar\'any \cite{Barany82}, we obtain colorful starshapedness of the $L_p$-sum, assuming we have at least $n+1$ colors.

\begin{restatable}{theorem}{starshapedness}
\label{thm:colorfulstarshapedness+prop11}
    Let $k\geq n+1$. Let $A_1,...,A_k\subset \R^n$ be closed sets such that $0\in \cap_{i=1}^k \textup{conv}(A_i)$. Then for any $p\in (1,\infty]$, the sum $A_1 \oplus_p \cdots \oplus_p A_k$ is starshaped.
\end{restatable}
In particular, starshapedness implies that $A^{\oplus_p k}\cap \R u$ is a segment, for any line $\R u$. This remains true even when dropping the assumption $0\in \text{conv}(A)$. 
\begin{restatable}{theorem}{lineintersectionisconvex}
\label{prop:lineintersectionisconvex}
Let $p\in (1,\infty)$ and  $A\subset \R^n$ be a compact set. Let $k\geq n+1$. Then for any unit vector $u\in \R^n$, $\R u\cap A^{\oplus_p k}$ is a segment. 
\end{restatable}
Note that this also holds true for $p=\infty$, since $A^{\oplus_{\infty} k}=\text{conv}(A)$ for $k\ge n+1$.

\subsection{Hausdorff metric}
For a compact $A\subset\mathbb{R}^n$, define
\begin{align*}
    d(A):=d_H(A,\textup{conv}(A))
\end{align*}
to be the distance between $A$ and its convex hull with respect to the Hausdorff metric $d_H$. The Shapley--Folkman--Starr theorem \cite[Appendix]{ST-69} implies that 
\begin{align}\label{eq:convergence_minkowski}
    \lim_{k\rightarrow\infty}d\left(\frac{1}{k}\sum_{i=1}^{k}A\right)=0,
\end{align}
 with convergence rate $1/k$, or in other words, that the sequence of averages $k^{-1}\sum_{i=1}^{k}A$ of a compact set $A$ converges to $\textup{conv}(A)$ in the Hausdorff metric with the same rate $1/k$.

 In the setting of Firey $L_p$-summation with $p>1$, in view of the following theorem, it is natural to consider the distance between $A$ and $\textup{conv}(0,A)$, which we denote by
\begin{align*}
    d_0(A):=d_H(A,\textup{conv}(0,A)).
\end{align*}

\begin{restatable}{theorem}{FireyConvergence}
\label{thm:Firey_convergence}
    Fix $1<p< \infty$. If $A\subset\mathbb{R}^n$ is a nonempty compact set, then
    \begin{align*}
        \lim_{k\rightarrow\infty}d_0\left(\frac{1}{k^{1/p}}A^{\oplus_pk}\right)=0,
    \end{align*}
    with convergence rate $k^{-1/p}$. The convergence is eventually monotone.
\end{restatable}
Notice that in the above theorem the appropriate normalization is $k^{-1/p}$ rather than $k^{-1}$. Indeed, this choice of normalization is further justified by observing that if $A$ is a convex set that does not necessarily contain the origin, then it is not difficult to see that for $1< p < \infty$, $k^{-1/p}A^{\oplus_pk}=[k^{-1/p},1]A\to [0,1]A=\textup{conv}(0,A)$ in Hausdorff distance when $k\to+\infty$.
In the case $p=\infty$, Carath\'eodory's theorem implies that $A^{\oplus_{\infty}k}=\textup{conv}(A)$ for $k\geq n+1$, so the convergence part of Theorem \ref{thm:Firey_convergence} holds, but with $d$ instead of $d_0$. 

The proof of Theorem \ref{thm:Firey_convergence} is asymptotic in the sense that the sumset $A^{\oplus_pk}$ has convenient properties when $k$ is sufficiently larger than $n$. Accordingly, a question not answered by our proof of Theorem \ref{thm:Firey_convergence}, is what sharp bound is possible if we take the sum of two sets. In the paper \cite{DF-2005}, Dyn and Farkhi posed this type of question by conjecturing that if $A,B\subset\mathbb{R}^n$ are nonempty compact sets, then
\begin{align*}
    d(A+B)^2\leq d(A)^2+d(B)^2.
\end{align*}
The Dyn--Farkhi conjecture is trivial in the case $n=1$, and was disproved by Fradelizi, Madiman, Marsiglietti, and Zvavitch \cite[Theorem 7.3]{FMMZ-18} when $n\geq 3$. Recently, van Hintum showed in \cite{VH-2026} that there are counterexamples even in the symmetric case $A=B$. The third-named author verified that the conjecture is true \cite{M-2025} when $n=2$, and later studied it in the setting of the $\ell_p$-norm \cite{M-2026} and arbitrary norms \cite{M2-2026}.

For the case $n=2$, we will prove the following bound for the Hausdorff distance of the $L_p$-sum of compact sets.

\begin{restatable}{theorem}{HausdorffDimTwoBound}\label{theorem:lp_hausdorff_dimension_2_max_bound}
    Fix $1<p\leq \infty$. If $A,B\subset\mathbb{R}^2$ are nonempty compact sets, then
    \begin{align*}
        d_0(A\oplus_p B)\leq \max\{d_0(A),d_0(B)\}.
    \end{align*}
    The bound is sharp in the sense that there exist nonconvex $A$ and $B$ for which equality holds.

    For $p=\infty$, we have the bound
    \begin{align*}
        d(A\oplus_{\infty}B)\leq \max\{d(A),d(B)\}.
    \end{align*}
    There exist nonconvex $A$ and $B$ for which equality holds.
\end{restatable}

In the case $n\geq3$, we have the following bound which generalizes \cite[Theorem 7.1]{FMMZ-18}.

\begin{theorem}\label{theorem:Lp_hausdorff_p_subadditive}
    Fix $1< p<\infty$. If $A,B\subset\mathbb{R}^n$ are nonempty compact sets, then 
    \begin{align*}
        d_0(A\oplus_pB)^p\leq d_0(A)^p+d_0(B)^p.
    \end{align*}
    In the case of $p=\infty$, we have
    \begin{align*}
        d_0(A\oplus_\infty B)\leq \max\{d_0(A),d_0(B)\}.
    \end{align*}
\end{theorem}

Moreover, in the above theorem, the power $p$ is optimal (see Proposition \ref{proposition:sharp_dimension_larger_than_3_hausdorff_lp_sum}).

\subsection{Effective standard deviation}\label{section:effective_standard_deviation}

Another measure of interest is the effective standard deviation, defined by Cassels in \cite{C-1975}, which gives a probabilistic measure of nonconvexity. For a compact $A\subset\mathbb{R}^n$, the effective standard deviation is defined by
\begin{align*}
    v^2(A)=\sup_{x\in \textup{conv}(A)}\inf\{V(X):\textup{supp}(X)\subset A,\mathbb{E}[X]=x\},
\end{align*}
where $V(X)$ is the variance of the random variable $X$. 
Equivalently, if we set
\begin{align*}
    v^2(x,A):=\inf\left\{\sum_{i=1}^{m} p_i\|x-a_i\|_2^2:x=\sum_{i=1}^{m} p_ia_i, p_i>0,\sum_{i=1}^{m} p_i=1,a_i\in A,m\in\mathbb{N}\right\},
\end{align*}
then we arrive at the equivalent definition
\begin{align*}
    v(A)=\sup_{x\in \textup{conv}(A)}v(x,A).
\end{align*}
Cassels proved that if $A,B\subset\mathbb{R}^n$ are nonempty compact sets, then
\begin{align}\label{eq:square_subadditive_cassels}
    v(A+B)^2\leq v(A)^2+v(B)^2.
\end{align}
By induction, \eqref{eq:square_subadditive_cassels} implies that the sequence of averages $k^{-1}\sum_{i=1}^{k}A$ converges to $\textup{conv}(A)$ in the measure $v$. We will prove the following $L_p$-sum analogue of \eqref{eq:square_subadditive_cassels}.
\begin{restatable}{theorem}{Pcassel}
\label{prop:pcassel}
    Let $1<p\leq \infty$ and let $n$ be a positive integer. Let $A,B\subset \R^n$ be nonempty compact sets such that $0\in \textup{conv}(A)\cap \textup{conv}(B)$. For $n\geq 3$, one has
    \begin{align*}
        v(A\oplus_p B)\leq 
        \begin{cases}
            (v(A)^{2p/(2-p)}+v(B)^{2p/(2-p)})^{(2-p)/2p}, & \textup{if } 1< p<2,\\
            \max\{v(A),v(B)\}, & \textup{if } 2\leq p\leq \infty.
        \end{cases}
    \end{align*}
    There exist nonconvex $A$ and $B$ so that these bounds are sharp. \\
    For $n=2$, one has the sharp bound $v(A\oplus_p B)\leq\max\{v(A),v(B)\}$ for every $p>1$.
\end{restatable}

It follows from Theorem \ref{prop:pcassel} that if $A$ is a nonempty compact set containing $0$, then the sequence of $L_p$-sum averages $k^{-1/p}A^{\oplus_pk}$ converges to $\textup{conv}(A)$ in the measure $v$. Furthermore, we show in Theorem \ref{proposition:strong_fractional_subadditivity_v} that the convergence is eventually monotone.

\subsection{Volume} Let $|A|$ denote the Lebesgue measure of a set $A\subset\mathbb{R}^n$. The volume deficit $\Delta(A):=|\textup{conv}(A)\backslash A|$ is a natural measure of nonconvexity. It is known from work of Emerson and Greenleaf \cite{EG-1969}, and Fradelizi, Madiman, Marsiglietti, and Zvavitch \cite[Theorem 3.7]{FMMZ-18}, that if $A$ is compact with nonempty interior, then the averages $k^{-1}\sum_{i=1}^{k}A$ converge to the convex hull of $A$ with respect to the volume deficit. We will prove the following analogous result for the $L_p$-sum.
\begin{restatable}{theorem}{ConvergenceVolume}
\label{thm: convergence-volume}
    Let $1<p< \infty$. If $A\subset \mathbb{R}^n$ is a compact set, then
    \begin{align*}
        \lim_{k\rightarrow\infty}\left|\frac{1}{k^{1/p}}A^{\oplus_pk}\right|=|\textup{conv}(0,A)|.
    \end{align*}
\end{restatable}

It is natural to ask whether the convergence is monotone. In the case of Minkowski sum, in the papers \cite{FMMZ-16,FMMZ-18}, Fradelizi, Madiman, Marsiglietti, and Zvavitch determined that the convergence is monotone when $n=1$, and not monotone when $n\geq 12$. It was later shown by Fradelizi, L\'angi, and Zvavitch \cite{FLZ-22} that monotonicity holds for starshaped sets when $n=2,3$, and that there is eventual monotonicity for starshaped sets when $n\geq 4$. The question of monotonicity remains open for general compact sets in the case $2\leq n\leq 11$. For the $L_p$-sum, we give the negative answer to the monotonicity of volume for any $n \geq 2\left\lceil\frac{p \ln 2}{\ln (9 / 8)}\right\rceil$. For $1<p<\infty$, we prove that the normalized volumes are eventually nondecreasing for every nonempty compact set.

\begin{restatable}{theorem}{MonotoneCol}
    Let $ n \geq 2$, $ 1 < p < \infty$, and let $A $ be a compact set. Then, for every $k \geq \max \{n,(n-1)(n-2)\}$,
    \begin{equation} \label{eq: monotone-col-compact}
        \left|\frac{A^{\oplus_p(k+1)}}{(k+1)^{1 / p}}\right| \geq\left|\frac{A^{\oplus_p k}}{k^{1 / p}}\right| .
    \end{equation}
    Equality holds if and only if $A$ lies in a proper linear subspace or $A^{\oplus_p k} =k^{1 / p} \operatorname{conv}(0,A) .$
\end{restatable}

\subsection{Related research} The work of Wegmann, covered in the survey \cite{FMMZ-18} and the paper \cite{W-1980}, establishes relations between some measures of nonconvexity. We also point out the paper by Schneider \cite{S-1975}, which defines another measure of nonconvexity through the Minkowski sum of a set with its convex hull. See also \cite{M3-2025,AFS-2025}. Other measures include the Beer index \cite{Beer-1973}, convexity ratio \cite{Goodman-1981,Balko-2017}, and reach \cite{Federer-1959}.

The convexification phenomenon has been studied in other settings. For example, \cite{AK-2025,Artstein-2021} focus on the setting of Banach spaces and infinite-dimensional spaces.

\subsection*{Acknowledgments}
The authors thank the LAMA at  University Gustave Eiffel and the functional analysis team of the Mathematical Institute of Jussieu for their hospitality.

\subsection*{Funding} 
M. Meyer is supported by the NSF through the MSPRF program (award number:
2502794). 
A.~Manui and A.~Zvavitch were partly supported by the U.S. National Science Foundation Grant DMS-2247771 and the United States--Israel Binational Science Foundation (BSF) Grant 2018115.
A.~Manui was supported by the Chateaubriand Fellowship of the Office for Science \& Technology of the Embassy of France in the United States.

\subsection*{AI disclosure statement}
The authors used an LLM as a research-assistance tool to suggest auxiliary lemmas, to check for misprints, and to improve the exposition of this article. Namely, an LLM provided significant assistance with the following content: it shortened or improved the proofs we had of Theorem \ref{prop:lineintersectionisconvex}, Theorem \ref{thm: convergence-volume}, and Lemma \ref{lem:containing-hyperplanebis}, it refined the proof of Theorem \ref{thm: hausd-convergence-2}, it provided the example of the equality case for $1<p<2$ in Theorem~\ref{prop:pcassel} and the example in Remark \ref{rem:Dis-Hau}. It also suggests Lemma \ref{lemma:firey_parallelogram_estimate}.


\section{Preliminaries}
Denote by $\{e_1,\ldots,e_n\}$  the standard basis of $\R^n.$ 
Recall the definition of $L_p$-sums of compact sets given in Lutwak, Yang, and Zhang \cite{LYZ-12}: for  $1< p \leq \infty$ and for any compact sets $K, L \subset \mathbb{R}^n$, define
\begin{equation}
\label{def:LYZ}
K\oplus_p L=\left\{(1-t)^{1 / q} x+t^{1 / q} y: x \in K, y \in L, \text { and } 0 \leqslant t \leqslant 1\right\},
\end{equation}
where $q$ denotes the H\"older conjugate of $p$ (i.e., $\frac{1}{p}+\frac{1}{q}=1$). For $p=1$, we have $K \oplus_1 L = K+L$. We note that the sets under consideration may not contain the origin. Consequently, the monotonicity properties via inclusion may fail in this setting: for $p_1 <p_2$, $K\oplus_{p_1} L$ is not necessarily contained in $K\oplus_{p_2} L$, as seen in the picture below.
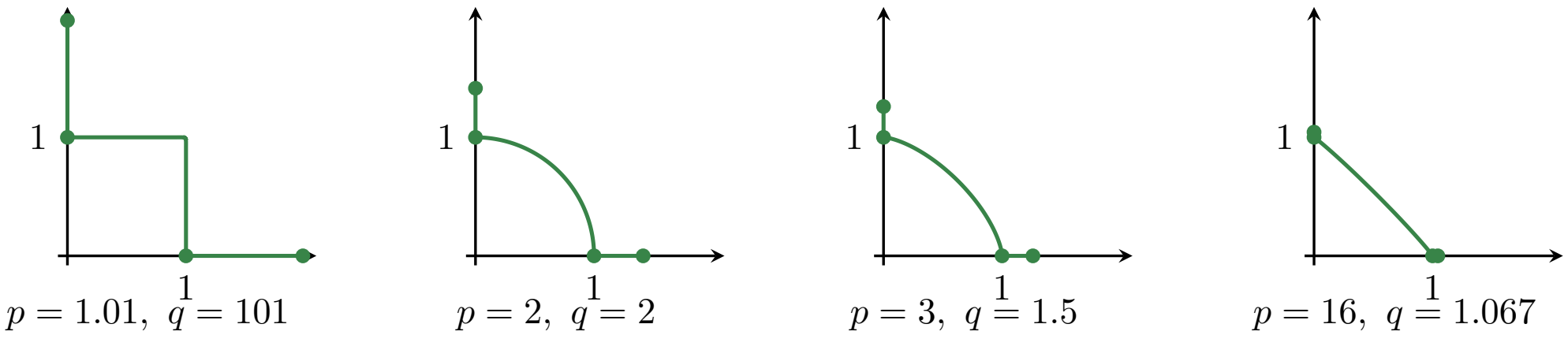
\begin{figure}[H]
\begin{tikzpicture}[scale = .8]
    \LqPicture{1.01}{0cm}
\end{tikzpicture}
\qquad 
\begin{tikzpicture}[scale = .8]
    \LqPicture{2}{0cm}
\end{tikzpicture}
\qquad 
\begin{tikzpicture}[scale = .8]
\LqPicture{3}{-6.2cm}
\end{tikzpicture}
\qquad 
\begin{tikzpicture}[scale = .8]
\LqPicture{16}{-12.4cm}
\end{tikzpicture}
\caption{$K\oplus_p L$ for $K=L= \{e_1,e_2\}$ and different values of $p$.}
\end{figure}
\noindent
The $L_p$-sum of compact sets need not be convex, as shown above. But notice that for $p >1$, we always have $K\cup L\subset K\oplus_p L$. 
Moreover, for $1<p\leq\infty$, the $L_p$-sum of any two nonempty sets $K,L\subset\mathbb R^n$ is path-connected. Indeed, for every $a\in K$ and $b\in L$, the continuous curve
\[
\gamma_{a,b}(t)=(1-t)^{1/q}a+t^{1/q}b, \qquad t\in[0,1],
\]
joins $a$ to $b$ within $K\oplus_p L$. Thus, any two points of $K\oplus_p L$ can be connected by joining at most three such paths. 
We note that the $L_p$-sum is associative and commutative, as follows from its definition. For convenience, we denote 
\[
    \underbrace{A_1 \oplus_p \cdots \oplus_p A_1}_{m_1 \textup{ copies}} \oplus_p \ldots \oplus_p \underbrace{A_k \oplus_p \cdots \oplus_p A_k}_{m_k \textup{ copies}}
    =
    \bigoplus^{(p)}_{j \in [k]} A_j^{\oplus_pm_j} .
\]

\begin{lemma}\label{lem:coefficient-reduction}
Let $A\subset\R^n$ be a nonempty compact set and $m = \dim \Span (A)$. 
Suppose that
\[
x=\sum_{j=1}^{r}\lambda_j a_j,
\qquad a_j\in A,\quad \lambda_j\geq0.
\]
Then there exists $\{i_1,...,i_s\}\subset [r]$ of size $s\le m$, such that $a_{i_1}, ...,a_{i_s}$ are linearly independent, and such that $x=\sum_{j=1}^s c_j a_{i_j}$ for some $c_j \ge 0$, with $||c||_1\le ||\lambda||_1$.

\end{lemma}

\begin{proof}
We may assume that $x\neq 0$, $\lambda_j>0$ for $1\le j \le r$, and that the $a_j$ are linearly dependent. So there exist real numbers
$\mu_1,\ldots,\mu_r$, not all zero, such that
$\sum_{j=1}^{r}\mu_j a_j=0,$ and  $\sum_{j=1}^{r}\mu_j\geq 0.$
In particular, at least one $\mu_j$ is positive. Let $I=\{t\in \R : \lambda+t\mu \in \mathbb R_+^r\}$. Then $I=[-t_0,t_1]$ for some real $t_0>0$ and with $t_1\in (0,+\infty]$. For $1\le j\le r$, set $c_j=\lambda_j-t_0\mu_j$. Then, $c_j\geq0$ for every $j$, and at least one $c_j$ vanishes. Moreover

\[
\sum_{j=1}^{r}c_j a_j
=\sum_{j=1}^{r}\lambda_j a_j
-t_0\sum_{j=1}^{r}\mu_j a_j=x,
\quad 
\sum_{j=1}^{r}c_j
=\sum_{j=1}^{r}\lambda_j
-t_0\sum_{j=1}^{r}\mu_j
\leq\sum_{j=1}^{r}\lambda_j.
\]
and so $|\text{supp}(c)|<|\text{supp}(\lambda)|$ and $||c||_1\leq ||\lambda||_1$. We repeat this argument until the remaining $a_i$'s are linearly independent.
\end{proof}
\begin{remark} In the above proof, if $0\notin \text{conv}(A)$ then $t_1<+\infty$. So in that case, we could choose $t_1$ rather than $t_0$ at each step, and then we end up with a representation of $x=\sum_{j\in J} c_j a_j$ with $J\subset [r]$ and the $(a_j, j\in J)$ linearly independent, and with $||c||_1\ge ||\lambda||_1$. The following lemma generalises this observation.
\end{remark}
\begin{lemma} \label{lem:coefficient-reduction-q}
    Let $A\subset\R^n$ be a nonempty compact set and $m = \dim \Span (A)$. Suppose that
\[
x=\sum_{j=1}^{r}\lambda_j a_j,
\qquad a_j\in A,\quad \lambda_j\geq0.
\]
If $ 0 \not \in \conv (A)$, then for any $q \in [1,\infty]$, there exist $\{i_1,...,i_s\}\subset [r]$ of size $s\le m$, such that $a_{i_1}, ...,a_{i_s}$ are linearly independent, and such that $x=\sum_{j=1}^s c_j a_{i_j}$ for some $c_j \ge 0$, with $||c||_q\ge ||\lambda||_q$.
\end{lemma}

\begin{proof}
   We may assume that $\lambda_j>0$ for $1 \leq j \leq r$, and that $a_j$ are linearly dependent. We can write $\sum_{i=1}^r \mu_i a_i =0$. Since $ 0 \not\in \conv(A)$, $\mu_1,\ldots, \mu_r$ have both positive and negative entries. Hence there exists $t_+, t_->0$ such that $I:=\{ t \in \R: \lambda + t\mu \in \R^r_+\}=[-t_-, t_+]$. 
   Let us write $c_j^{\pm}:=\lambda_j \pm t_{\pm} \mu_j$. This yields two distinct representations $x=\sum_{j=1}^r c_j^{\delta} a_j$, where $\delta\in\{\pm\}$, with $c_j^{\delta}\ge 0$ and with $|\text{supp}(c^{\pm})|\leq r-1$. Moreover, by convexity of the $\ell_q$-norm, 
   since $0\in \text{int}(I)$, we can pick $c\in \{c^+,c^-\}$ with $||c||_q \ge  ||\lambda||_q$.
   We can repeat this step-by-step elimination until the remaining $a_i$'s are linearly independent.
 \end{proof}

We prove an $L_p$-sum version of the Shapley--Folkman lemma.
\begin{lemma}[The Shapley--Folkman lemma for $\oplus_p$]\label{shapley_folkman_LP}
Let $1 \leq p \leq \infty$, $k\geq n+1$ and let $A_1, \ldots, A_k \subset \R^n$ be nonempty sets not necessarily containing the origin. If $x\in \bigoplus_{i\in [k]}^{(p)} \textup{conv}(A_i)$, then there exists $I\subset [k]$, of cardinality at most $n$, such that 
\[x\in \left(\bigoplus^{(p)}_{i\in I} \textup{conv}(A_i)\right) \oplus_p \left(\bigoplus^{(p)}_{i\in I^c} A_i \right).\]
\end{lemma}
\begin{proof}
The case $p=1$ is the classical Shapley--Folkman lemma. Thus, we may assume that $1<p\leq\infty$. A point $x \in \bigoplus_{i\in [k]}^{(p)} \operatorname{conv}(A_i)$ can be written as
\[
x=\sum_{i=1}^k t_i^{1/q}a_i,
\qquad
a_i\in\operatorname{conv}(A_i),\quad
t_i\geq 0,\quad \sum_{i=1}^k t_i=1.
\]
Then, $x\in\bigoplus_{i\in[k]}^{(1)} \operatorname{conv}(t_i^{1/q}A_i).$ By the classical Shapley--Folkman lemma, there exists
$I\subset[k]$, with $|I|\leq n$, such that
\[
x\in
\left(\bigoplus_{i\in I}^{(1)}
\operatorname{conv}(t_i^{1/q}A_i)\right)
+
\left(\bigoplus_{i\in I^c}^{(1)}
t_i^{1/q}A_i\right)
\subset
\left(\bigoplus_{i\in I}^{(p)}
\operatorname{conv}(A_i)\right)
\oplus_p
\left(\bigoplus_{i\in I^c}^{(p)}A_i\right),
\]
where the last inclusion follows from the definition of $L_p$-sums. 
\end{proof}

\section{Geometry of \texorpdfstring{$L_p$}{Lp}-sums}\label{section:goemetry_of_lp_sums}

\subsection{Properties of the convex hull of the \texorpdfstring{$L_p$}{Lp}-sum}
If $A,B\subset\mathbb{R}^n$ are nonempty sets, then the identity $\textup{conv}(A+B)=\textup{conv}(A)+\textup{conv}(B)$ is useful for comparing a sumset to its convex hull, as observed in the paper \cite{FMMZ-18}. The next result gives the $L_p$-sum analogue of this identity.
\begin{lemma}
\label{lemma:conv0p}
    Let $n\ge 1$ be an integer. Let $A,B$ be compact sets in $\mathbb{R}^n$. Then, if $1<p<\infty$:
    \[
    \textup{conv}(0,A\oplus_p B)=\textup{conv}(0,A)\oplus_p\textup{conv}(0,B).
    \] 
    and for $p=\infty$:
    \begin{equation}
    \label{eqn:convinfty}
    \textup{conv}(A\oplus_{\infty}B)=\textup{conv}(A, B)=\textup{conv}(A)\oplus_{\infty}\textup{conv}(B).
    \end{equation}
\end{lemma}

\begin{proof}
    As the $L_p$-sum of two convex sets containing $0$ is convex, the inclusion $\subset$ is clear, so let us justify the other one. Denote $K=\text{conv}(0,A)$ and $L=\text{conv}(0,B)$.  Let $z\in K\oplus_p L$. We can assume $z\neq 0$. We can write $z=\alpha x+\beta y,$ for some $x\in K, y\in L,$ $\alpha,\beta\geq0,$ with $\alpha^q+\beta^q=1.$
There exist $s,t\in[0,1]$ and convex combinations
\[
x=s\sum_{i=1}^{m_A}\lambda_i a_i,
\qquad
y=t\sum_{j=1}^{m_B}\mu_j b_j,
\]
where $a_i\in A$, $b_j\in B$, and
$\lambda_i,\mu_j\geq0$, with
$\sum_i\lambda_i=\sum_j\mu_j=1$.
Set
\[
r=\big((\alpha s)^q+(\beta t)^q\big)^{1/q}\leq1 \, \quad \widehat\alpha=\frac{\alpha s}{r}, \quad 
\widehat\beta=\frac{\beta t}{r}
\]
which are well-defined since $z\neq 0$ implies that $r>0$.
Then $\widehat\alpha^q+\widehat\beta^q=1$, and
\[
\begin{aligned}
z
&=r\left(
  \widehat\alpha\sum_i\lambda_i a_i
  +\widehat\beta\sum_j\mu_j b_j
  \right)=r\sum_{i,j}\lambda_i\mu_j
  \left(\widehat\alpha a_i+
        \widehat\beta b_j\right).
\end{aligned}
\]
Since each
$\widehat\alpha a_i+\widehat\beta b_j$
belongs to $A\oplus_p B$, and since
$\lambda_i\mu_j\geq0,$
$\sum_{i,j}\lambda_i\mu_j=1,$ $r\in [0,1]$,
we deduce that $z\in\text{conv}(0, A\oplus_p B)$.
This proves the reverse inclusion.

We now prove the equalities \eqref{eqn:convinfty}. Clearly, 
$A\cup B \subset A\oplus_{\infty} B\subset \text{conv}(A, B).$
Taking convex hulls gives $\text{conv}(A, B)=\text{conv}(A\oplus_{\infty} B)$. The inclusion $\text{conv}(A)\oplus_{\infty} \text{conv}(B) \subset \text{conv}(A, B)$ is clear, so we only need to justify the other inclusion. Let $x\in \text{conv}(A,B)$. We can write $x$ as a finite convex combination:
\[
x=\sum_i\alpha_i a_i+\sum_j\beta_j b_j,
\qquad
\alpha_i,\beta_j\ge0,
\qquad
\sum_i\alpha_i+\sum_j\beta_j=1,
\]
with $a_i\in A$ and $b_j\in B$. Set $t=\sum \alpha_i$, so $1-t=\sum \beta_i$. If $t=0$, then $x\in \text{conv}(B)$; if $t=1$, then $x\in \text{conv}(A)$. If $0<t<1$, then
\[
x=
t\left(\sum_i\frac{\alpha_i}{t}a_i\right)
+
(1-t)\left(\sum_j\frac{\beta_j}{1-t}b_j\right) =: t a' +(1-t) b' \in \text{conv}(A) \oplus_{\infty} \text{conv}(B)
\]
showing the other inclusion.
\end{proof}
    
Notice that the above equality does not hold for $p=1$: take $A=\{e_1\}$ and $B=\{e_2\}$.


\subsection{Compact sets}



A simple example of a nonconvex compact set that one can study is the set $A=\{e_1,\dots,e_n\}$. The iterated $L_p$-sum of $A$ with itself can be computed exactly.  
\begin{lemma} \label{lem:sum-of-ei}
    Fix \(1< p< \infty\). For every integer \(k\geq 1\), we have
    \begin{align}
        \{e_1,\ldots,e_n\}^{\oplus_p k} = 
        \bigcup_{\substack{m_1,\ldots,m_n\in \mathbb{Z}_{\geq 0} \\ \sum_{i=1}^n m_i =k} }\left(\diag (m_1^{1/p},\ldots,m_n^{1/p}) (B_q^n \cap \R^n_+)\right) \setminus \left(\Int (B_q^n )\cap \R^n_+\right).
    \end{align}
\end{lemma}

\begin{proof}
    Using \eqref{def:LYZ}, we have 
    \begin{align}
        \{e_1,\ldots,e_n\}^{\oplus_p k} &= \bigcup_{\substack{(m_1,\ldots,m_n) \\ \sum_{i=1}^n m_i =k} } \,
        \bigoplus_{i \in [n]}^{(p)}  \{e_i\}^{\oplus_p m_i}.
    \end{align}
    Thus, it is enough to prove that for any nonnegative integers $m_1,\ldots,m_n$ such that $\sum_{i=1}^n m_i = k$, we have
    \[
        \bigoplus_{i \in [n]}^{(p)}  \{e_i\}^{\oplus_p m_i}
        = 
        \left(\diag (m_1^{1/p},\ldots,m_n^{1/p}) (B_q^n \cap \R^n_+)\right) \setminus \left(\Int (B_q^n )\cap \R^n_+\right).
    \]
    If some of the $m_i$'s are zero, the same argument applies in the corresponding lower-dimensional subspace. We then embed the result back into the original space by inserting zeros in the coordinates corresponding to those indices.
    
    Assume that $m_i \geq 1$ for all $i.$ First, for each $i$, we have
 \[
        E_i:=
       \{e_i\}^{\oplus_p m_i}
        = \left\{\sum_{j=1}^{m_i} t_j^{1 / q} e_i: 0 \leqslant t_j \leqslant 1 \text { and } \sum_{j=1}^{m_i} t_j = 1\right\} = [e_i, m_i^{1/p}e_i],
    \]
    by H\"older's inequality and by the path-connectedness of $L_p$-sums.
    Let
    \[
        E:= \bigoplus_{j\in [n]}^{(p)} E_j  = \left\{\sum_{i=1}^{n} t_i^{1 / q} \lambda_i e_i: 1 \leq \lambda_i \leq m_i^{1/p}, 0 \leqslant t_i \leqslant 1 \text { and } \sum_{i=1}^{n} t_i = 1\right\}.
    \]
    We will prove that $E = \left(\diag \left(m_1^{1/p},\ldots,m_n^{1/p}\right)  (B_q^n \cap \R^n_+)\right) \setminus \left(\Int (B_q^n )\cap \R^n_+\right).$
    Let $x = \sum_{i=1}^{n} t_i^{1 / q} \lambda_i e_i$ be a point in $E$ where $t_i \in [0,1]$ with $\sum_{i=1}^n t_i = 1$, and $\lambda_i \in [1,m_i^{1/p}]$. Then,
    \[
        \left\| \sum_{i=1}^{n} \frac{t_i^{1 / q} \lambda_i}{m_i^{1/p}} e_i \right\|_q^q = \sum_{i=1}^{n} t_i \left(\frac{\lambda_i}{m_i^{1/p}}\right)^q \leq \sum_{i=1}^{n} t_i =1 .
    \]
    Using the fact that all points in $E$ are the sum of points in the positive orthant, we have $x \in \diag \left(m_1^{1/p},\ldots,m_n^{1/p}\right) (B_q^n \cap \R^n_+)$. Since 
    \[
        \left\| \sum_{i=1}^{n} t_i^{1 / q} \lambda_i e_i \right\|_q^q = \sum_{i=1}^{n} t_i \lambda_i^q \geq \sum_{i=1}^{n} t_i  = 1,
    \]
    we get $E \subset \left(\diag \left(m_1^{1/p},\ldots,m_n^{1/p}\right)  (B_q^n \cap \R^n_+)\right) \setminus \left(\Int (B_q^n )\cap \R^n_+\right).$
    Reciprocally, let $$x = (x_1,\ldots,x_n) \in \left(\diag \left(m_1^{1/p},\ldots,m_n^{1/p}\right)  (B_q^n \cap \R^n_+)\right) \setminus \left(\Int (B_q^n )\cap \R^n_+\right).$$ 
    We consider the function $f$ on $[0,1]$ defined by
    \[
        f(s) = \sum _{i=1}^n \left(\frac{x_i}{m_i^{s/p}}\right)^q.
    \]
    Since $f$ is continuous with $f(0) = \sum _{i=1}^n x_i^q \geq1$ and $f(1) = \sum _{i=1}^n \left(\frac{x_i}{m_i^{1/p}}\right)^q \leq 1$, by the intermediate value theorem, there exists $s_0 \in [0,1]$ such that $f(s_0)=1$.
    Let $ t_i = \left(\frac{x_i}{m_i^{s_0/p}}\right)^q$. Hence $\sum_{i=1}^n t_i =f(s_0)=1$ and 
    $ \sum_{i=1}^n t_i^{1/q} m_i^{s_0/p}e_i = \sum_{i=1}^n x_i e_i = x.$
    Set $\lambda_i = m_i^{s_0/p} $. Since $ 1 \leq \lambda_i \leq m_i^{1/p}$, we get $ x \in E$ and the proof is complete.
\end{proof}


\begin{prop}
    Fix $1< p < \infty$. For every compact set $A \subset \R^n$ such that $\Span(A)=\R^n$, we have 
    \[
    \Int\left(A^{\oplus_p(n+1)}\right) \neq \emptyset.
    \]
\end{prop}

\begin{proof}
   The set $A$ contains $n$ linearly independent vectors. By applying a linear transformation, we assume that $\{e_1,\ldots,e_n\} \subset A. $ 
    Using Lemma \ref{lem:sum-of-ei}, we have
    \begin{align}
       \left(\diag (1,\ldots,1,2^{1/p}) (B_q^n \cap \R^n_+)\right) \setminus \left(\Int (B_q^n )\cap \R^n_+\right) 
       &\subset \{e_1,\ldots,e_n\}^{\oplus_p (n+1)}
       \subset A^{\oplus_p (n+1)}
       .
    \end{align}
    The left set has nonempty interior, and so $A^{\oplus_p (n+1)}$ has nonempty interior.
    \end{proof}
    \begin{remark} \label{remark 1}
        We note that  $n+1$ is optimal: for $A=\{e_1,\ldots,e_n\}$ one has $\Int\left(A^{\oplus_p n}\right) = \emptyset.$ Indeed, it follows from Lemma~\ref{lem:sum-of-ei} that 
        $A^{\oplus_p n} \subset \left(\bigcup_{k=1}^n e_k^{\perp}\right) \cup \left((\partial B_q^n)\cap \R_+^n\right)$,
        and so $A^{\oplus_p n}$ has empty interior.
        
    \end{remark}

\begin{lemma} \label{lem:containing-hyperplanebis}
    Fix $1<p \leq \infty$. For every integer $k\ge n$, we have
\begin{equation}
\label{eqn:inclusionhyperplane}
        \left\{(x_1,\ldots,x_n)\in \R_+^n: \|x\|_q \geq 1, \|x\|_1 \leq (k-n+1)^{1/p} \right\} \subset  \{e_1,\ldots,e_n\}^{\oplus_p k}.
\end{equation}
Equality holds when $ p =\infty$. In particular, if $k\geq 2n-1$, we have
$$\left\{x\in \R_+^n: \|x\|_1 = (k-n+1)^{1/p} \right\} \subset  \{e_1,\ldots,e_n\}^{\oplus_p k}.$$
\end{lemma}
\begin{proof}
For the case $p=\infty$,  the assertion follows from Carath\'eodory's theorem. Thus, we may assume that $p < \infty$.
Let $ x =(x_1,\ldots,x_n) \in \R^n_+$ be such that $\|x\|_q \geq 1, $ and $ \sum_{i=1}^n x_i \leq (k-n+1)^{1/p}$.  For each $1\leq i \leq n-1$, let 
    \(
        m_i =  \lceil x_i (k-n+1)^{1/q} \rceil,
    \)
    and $m_n = k- \sum_{i=1}^{n-1} m_i$. Observe that
    \begin{align*}
    m_n=k-\sum_{i<n} m_i & \ge \, k-n+1 - \sum_{i<n} x_i (k-n+1)^{1/q} \geq  x_n (k-n+1)^{1/q}.
    \end{align*}
   It follows from the definition of $m_i$ that $ m_i = 0 $ only if $x_ i =0$.
    Therefore,
    \[
        \sum_{i: m_i >0} \left( \frac{x_i}{m_i^{1/p}} \right)^q \leq \sum_{i=1}^n \frac{x_i}{(k-n+1)^{1/p}} \leq  1.
    \]
    Thus, 
    $$
        x \in \diag (m_1^{1/p},\ldots,m_n^{1/p}) (B_q^n \cap \R^n_+).
    $$ 
    Since $\|x\|_q \geq 1 $, $x \notin \Int (B_q^n )\cap \R^n_+$. Using Lemma~\ref{lem:sum-of-ei}, we get that 
    \[
        x \in \left(\diag (m_1^{1/p},\ldots,m_n^{1/p}) (B_q^n \cap \R^n_+)\right) \setminus \left(\Int (B_q^n )\cap \R^n_+\right) \subset \{e_1,\ldots,e_n\}^{\oplus_p k}.  
    \]
If $k\ge 2n-1$ then $(k-n+1)^{1/p}\ge n^{1/p}$. Thus if $x\in \R_+^n$ and $||x||_1=(k-n+1)^{1/p}$, then $n^{1/p} \le ||x||_1\le n^{1/p} ||x||_q$, by H\"older's inequality. Thus, the stated inclusion for that case follows from the inclusion \eqref{eqn:inclusionhyperplane}.
\end{proof}

\subsection{Starshapedness results}
\label{section:segments}

The aim of this section is to show two results, Theorem~\ref{thm:colorfulstarshapedness+prop11} and Theorem \ref{prop:lineintersectionisconvex}. Both assume $k\geq n+1$. The former says that when $0\in \cap \text{conv}(A_i)$, then (for any $p>1$) $A_1\oplus_ p \cdots \oplus_p A_k$ is starshaped. In particular, this implies that $\R u \cap A^{\oplus_p k}$ is always an interval if $0\in \text{conv}(A)$. Theorem \ref{prop:lineintersectionisconvex} says that this weaker conclusion holds without any assumption on $A$.

We start by quoting two classical results from combinatorial geometry, due to B\'ar\'any~\cite{Barany82}, which are generalizations of Carath\'eodory's theorem (and of its cone version). They are known as \emph{colorful Carath\'eodory's theorems}, see \cite{Matousek}, Chapter 8.



\begin{fact}[B\'ar\'any {\cite[Theorem~2.1]{Barany82}}]
\label{fact:rainbowsimplices}
Let $k\geq n+1$ and $A_1, ... , A_k \subset \R^n$. Then $\cap_{i=1}^k \text{conv}(A_i) \subset \bigoplus_{i \in [k]}^{(\infty)} A_i =\{\sum_{i=1}^k \lambda_i a_i : a_i\in A_i, \lambda_i \geq 0, \sum \lambda_i=1\}$.
\end{fact}
The following is the cone-analogue of Fact \ref{fact:rainbowsimplices}, it appears as Theorem 2.2 in \cite{Barany82}.

\begin{fact}[B\'ar\'any {\cite[Theorem~2.2]{Barany82}}]
\label{fact:colorfulcones}
Let $k'\geq m$ and $A_1, ... , A_{k'} \subset \R^m$. Assume $v\in \R^m\setminus \{0\}$ is such that $v\in \cap_{i\le k'}\text{cone}(A_i)$. Then there exists $(a_i\in A_i,\mu_i\geq 0 , i\in [k'])$  such that $v= \sum_{i=1}^{k'} \mu_i a_i$. More precisely, there exists $1\leq r\leq m$, $J=\{c_1,...,c_r\}\subset [k']$, $a_1,...,a_r$ linearly independent, with $a_i \in A_{c_i}$, and $(\mu_i>0, i\in [r])$, such that $v=\sum_{i=1}^r \mu_i a_i$.
\end{fact}

We refer to \cite{Barany82} and to \cite[Theorem~8.2.1, p.~199]{Matousek} for proofs of the above two theorems. The latter introduces the notion of rainbow simplices, which we recall now, as we will need it for the next Proposition.
A $d$-simplex $T=\text{conv}(a_1,...,a_{d+1})$ is called rainbow, if all $d+1$ vertices have distinct colors. Since a vertex may have several colors, we make this definition more precise by calling $T=\text{conv}(a_1,...,a_{r})$ a $J$-rainbow simplex, if $\dim(T)=r-1$ and if there exists a bijective $\sigma: [r] \to J$ such that $a_i\in A_{\sigma(i)}$ for all $i\in [r]$. We denote by $\mathcal{T}_J$ the set of $J$-rainbow simplices, and we denote \[\mathcal{T}_{\ell}=\bigcup_{\substack{J\subset [k] \\ |J|=\ell+1}}\mathcal{T}_J\]
the set of $\ell$-dimensional rainbow simplices with colors in $[k]$.
Finally, denote by $\mathcal{T}=\cup_{\ell \ge 0} \mathcal{T}_{\ell}$ the set of all rainbow simplices with vertices in $\cup_{j\in [k]} A_j$.

Similarly, if $C=\text{cone}(a_1,...,a_r)=\R_+ a_1+\cdots+\R_+ a_r$ is an $r$-dimensional cone with the $a_i$'s \emph{admitting distinct colors}\footnote{i.e., as before, there exists $J\subset [k]$ of size $r$, and $\sigma : [r]\to J$ a bijection, such that $a_i\in A_{\sigma(i)}$ $\forall i\in [r]$.}, then we call $C$ an $r$-rainbow cone.

The following is the case $p=\infty$ of Theorem \ref{thm:colorfulstarshapedness+prop11}.
\begin{prop}
\label{prop11}
Let $k \geq n+1$ and let $A_1, \ldots, A_k\subset \R^n$ be arbitrary sets such that $x_0\in \cap_{i=1}^k \text{conv}(A_i)$. Then $R:= \bigoplus_{i \in [k]}^{(\infty)} A_i=A_1 \oplus_{\infty} \cdots \oplus_{\infty} A_k$ is starshaped at $x_0$.
\end{prop}
\begin{proof}
If we replace each $A_i$ with $A'_i:=A_i-x_0$, then both $\cap \text{conv}(A_i)$ and $A_1\oplus_{\infty} \cdots \oplus_{\infty} A_k$ get translated by $-x_0$, so we can assume $x_0=0$. By Carath\'eodory's theorem, we can assume that all $A_i$ are finite, and so $R:= A_1\oplus_{\infty} \cdots \oplus_{\infty} A_k=\cup_{T\in \mathcal{T}} T$ is compact, as a finite union of (compact) simplices. 
If $z\in R$, then consider $I(z):=\{t\in [0,1] : [tz,z]\subset R\}$. It is non-empty (as $z\in R$) and closed, since $R$ is compact.
Thus, by connectedness of $[0,1]$, it suffices to show that $I(z)$ is open (in $[0,1]$), to deduce $I(z)=[0,1]$ as needed. Therefore it suffices to show that if $x\in R\setminus \{0\}$, then there exists $\epsilon>0$ such that $[(1-\epsilon)x,x]\subset R$.


Fix $x\in R$, $x\neq 0$, and $a_1, \ldots , a_r$ with $r\geq 1$ minimal, such that $x\in T=\text{conv}(a_1,\ldots,a_r)$ and $T$ is rainbow.
Let us first look at the case when $r=1$, i.e. $x\in \cup_{i=1}^k A_i$, say $x\in A_1$, $x\neq 0$. There remain $k-1\ge n$ colors unused. By \cite[Theorem 2.3]{Barany82}, we can find $b_2\in A_2,..., b_{n+1}\in A_{n+1}$, such that $0\in \text{conv}(a_1,b_2,...,b_{n+1})$. Hence $[0,x]=[0,a_1]\subset T'$, with $T'=\text{conv}(a_1,b_2,...,b_{n+1}) \subset A_1\oplus_{\infty} A_2 \oplus_{\infty} \cdots \oplus_{\infty} A_{n+1} \subset R$.
%

Now assume $x\notin \cup_{i\in k} A_i$, so $r\geq 2$. Up to renaming colors, assume $x\in T=\text{conv}(a_1,\ldots ,a_r)$ with $a_i\in A_i$ (and $r$ minimal), and $T$ an $(r-1)$-rainbow simplex. By Carath\'eodory's theorem, $r\leq n+1$. If $r=n+1$, then minimality of $r$ means that $x\in \text{int}(T)$, in which case it is clear that $[(1-\epsilon)x,x]\subset R$ for some $\epsilon>0$. Hence assume $2\leq r\leq n$. So there remains $k-r\geq (n+1)-n=1$ colors. By minimality of $r$, we have $x\in \text{relint}(T)$. Let $W=\text{Span}(a_i-x ; i\leq r)=\text{Span}(T-x)$.  Denote $m=\text{dim}(W^{\perp})=n-r+1 \geq 1$. Let $J_r=[k]\setminus [r]$ be the colors not used by $T$.
 Denote $\pi$ the orthogonal projection onto $W^{\perp}\approx \R^m$.
For each $i\in J_r$, $\pi(-x)\in \text{cone}(\pi(A_i-x))$. Since $|J_r|=k-r\geq n+1-r=m$, Fact \ref{fact:colorfulcones} tells us that we can pick $J'\subset J_r$ of size at most $m$ (take $J'=\emptyset$ if $\pi x=0$) and
$(a_i, i\in J')$, $\mu_i > 0$, such that $-\pi x=\sum_{i\in J'} \mu_i \pi (a_i -x)$ and $\dim(\text{Span}(\pi(a_i-x) : i\in J'))=|J'|$. 
Write $-x=-\pi x+w_x$, and $a_i -x=\pi(a_i-x)+w_i$, with $w_x, w_i \in W$. Then, setting $w:=w_x-\sum_{i\in J'} \mu_i w_i \in W$, we get:
\begin{equation}
\label{eqn:stayinT}
-x=-\pi x+w_x=\sum \mu_i (a_i-x-w_i)+w_x=\sum_{i\in J'} \mu_i(a_i-x) +w.
\end{equation}
Set $T'=\text{conv}(T,(a_j , j\in J'))$. By linear independence of $(\pi(a_i-x)), i\in J')$, we see that
\begin{align*}
\dim(T')=\dim(\text{Span}(T'-x))&=\dim\left(\text{span}(a_i-x ; i\in [r])+\text{span}(a_j-x ; j\in J')\right)\\
&=r-1+\dim(\text{Span}(\pi(a_j-x) ;j\in J'))=r+|J'|-1\, .
\end{align*}
 In other words, $T'\in \mathcal{T}_{[r]\cup J'}$ is an $(r+|J'|-1)$-rainbow simplex. Set $M=\sum \mu_i >0$. For $\epsilon>0$ small enough, setting $\tilde{\mu_i}=\mu_i/M$ and $x'=x+\frac{\epsilon}{1-\epsilon M} w$, the identity \eqref{eqn:stayinT} gives :
\begin{equation}
\label{eqn:cancontinuetowards0}
(1-\epsilon) x=(1-\epsilon M)x+\epsilon w +\epsilon M \sum_{i\in J'} \tilde{\mu}_i a_i=: (1-\epsilon M)x' +\epsilon M \sum_{i\in J'} \tilde{\mu}_i a_i \quad.
\end{equation}
 If $\epsilon>0$ small enough, then $1-\epsilon M\in(\frac{1}{2},1)$ and $x+2\epsilon w\in T$ (since $w\in W=\text{Span}(T-x)$ and $x\in \text{relint}(T)$) so $x'\in [x,x+2\epsilon w] \subset T \subset T'$, so that \eqref{eqn:cancontinuetowards0} implies that $[(1-\epsilon)x,x]\subset T' \subset R$. We conclude that for any $x\in R$, there exists $\epsilon>0$ such that $[(1-\epsilon)x,x]\subset R$, which, by the connectedness argument mentioned earlier, implies that $R$ is starshaped.
\end{proof}

\begin{lemma}
\label{lem:easylemma}
   Let $r\ge 1$, $a_1,...,a_r\in \R^n$ and let $T=\text{conv}(a_1,...,a_r)$ and, for $p\in (1,\infty)$, let $T_p=  \{a_1\} \oplus_p \cdots \oplus_p \{a_r\}
   $. If $0\in T$, then $T\subset T_p$.
\end{lemma}
\begin{proof}
If $r=1$, then $0\in T$ means $a_1=0$, thus $T=T_p=\{0\}$. If $r=2$, then $T$ is a segment and one can check that $T=T_p$ in this case too. Assume $r\ge 3$. Then, since $0\in T$, we may write $0=\sum_{i=1}^r \eta_i a_i$ for some $\eta_i\geq 0$, $||\eta||_1=1$. Let $x\in T$. If $x=a_i$, then clearly $x\in T_p$. Assume $x=\sum_{i=1}^r \lambda_i a_i$, with $\lambda_i\ge 0$ and $||\lambda||_1=1$. We may assume $x\notin \{a_1,...,a_r\}$. Thus, $|\text{supp}(\lambda)|\ge 2$, and so $||\lambda||_q<1$. For any $t\ge 0$, write $x=\sum_{i=1}^r (\lambda_i +t\eta_i) a_i.$
Denote $\mu_i(t)=\lambda_i+t\eta_i \ge 0$, and so $||\mu(t)||_q=\left(\sum_{i=1}^r (\lambda_i +t\eta_i)^q\right)^{1/q}$ is increasing in $t$, starts with $||\mu(0)||_q=||\lambda||_q<1$, and goes to $\lim_{t\to +\infty} ||\mu(t)||_q=+\infty$. Hence by continuity, there exists $t_0 >0$ such that $||\mu(t_0)||_q=1$, and so $x=\sum_{i=1}^r \mu_i(t_0) a_i \in T_p$.
\end{proof}
\starshapedness*

\begin{proof} The case $p=+\infty$ follows from Proposition \ref{prop11} so we assume $1<p<+\infty$.
Denote by $q>1$ the conjugate exponent of $p$. 
If $A_1,...,A_k$ are finite sets, denote $R_p=A_1\oplus_p ... \oplus_p A_k$, which is a finite union of $T_p=T_p(a_1,...,a_\ell)=\{a_1\}\oplus_p ... \oplus_p \{a_\ell\}$, $\ell \leq k$, and similarly for $p=\infty$, $R_{\infty}$ is a finite union of $T=\text{conv}(a_1,...,a_\ell)$, $\ell \leq k$ (which, by Carath\'eodory's theorem, is the same as the finite union of those $T=\text{conv}(a_1,...,a_\ell)$ for which $T$ is an $(\ell-1)$-rainbow simplex). In Fact \ref{fact:rainbowsimplices}, to show that the assumptions $0\in \cap \text{conv}(A_i)$, and $k\geq n+1$, yield $0\in R_{\infty}$, one examines the quantity $\rho=\text{dist}(0,R_{\infty})=\min_{T\in \mathcal{T}} \text{dist}(0,T)$ with $\mathcal{T}=\cup_{0\le \ell \le n} \mathcal{T}_{\ell}$ the set of all rainbow simplices, and argues by a contradiction argument, that one cannot have $\rho>0$ (see, e.g., \cite[Section~8.2, proof of Theorem~8.2.1]{Matousek}).
The proof here is similar: if $z\in R_p^{*} :=\{\sum_{i\leq k} \lambda_i a_i : \lambda_i\geq 0, a_i\in A_i, \sum \lambda_i^q \leq 1\}$, we shall show that $z\in R_p$ by showing that $\max_{(a,\lambda)} F(\lambda)=1$, where $F(\lambda)=\sum_{i\in [k]} \lambda_i^q$ and the maximum is over all valid representations of $z\in R_p^*$, i.e. over all tuples $((a_i)_{i\in [k]}, (\lambda_i)_{i\in [k]})$ such that $\lambda_i\geq 0$, $a_i \in A_i$, $F(\lambda)\leq 1$ and $z=\sum \lambda_i a_i$. By Fact \ref{fact:rainbowsimplices}, there exists $T=T(a_1,...,a_r)$ a rainbow simplex (with $r\leq n+1$) such that $0\in T$, and hence (by Lemma \ref{lem:easylemma}) such that $0\in T_p=\{a_1\}\oplus_p \cdots \oplus_p \{a_r\}$, and so $0\in R_p$. So we can assume $z\neq 0$.

For contradiction, assume we had $m(z)<1$, where \[m(z):=\max_{a,\lambda} F(\lambda)=\max\left\{\sum_{i=1}^k \lambda_i^q : \lambda_i \geq 0, F(\lambda)\leq 1 ; \exists a_i\in A_i, z=\sum \lambda_i a_i \right\} \, ,\]
and let us fix $J\subset [k]$, $a_i\in A_i, \lambda_i >0$ (for $i\in J$), such that $(a,\lambda)$ is a valid representation of $z$, and $F(\lambda)=m(z)<1$ is maximal. Say $|J|=r\geq 1$ ($r>0$ because $z\neq 0$), and rename colours so that $J=[r]$. Firstly we claim that $a_1,...,a_r$ are linearly independent; in particular $r\leq n$.
Indeed, if we had $\sum h_i a_i=0$ for some $h\in \R^r$, $h\neq 0$, this would contradict maximality of $F(\lambda)$ : one has $0<\lambda_i<1$ for each $i\in [r]$, so $0<\lambda_i \pm \epsilon h_i <1$ (for all $i$)  and $F(\lambda\pm \epsilon h)<1$ for $\epsilon >0$ small enough (by continuity of $F$); in other words, for small $\epsilon>0$, $(a,\lambda+\epsilon h)$ and $(a,\lambda-\epsilon h)$ are also valid representations of $z$. By strict convexity of $F$, one has $\max\{F(\lambda+\epsilon h),F(\lambda-\epsilon h)\} > F(\lambda)$, a contradiction.

Assume we have $r=n$, i.e. $\R^n=\text{Span}(a_1,...,a_n)$. Hence $L(\sum h_i a_i):=\sum h_i \lambda_i^{q-1}$ defines a linear form on $\R^n$, and $L(z)=F(\lambda)>0$. Since $k\geq n+1$, we have at least one unused color. Since $0\in \text{conv}(A_{n+1})$, we can fix some $b=b_{n+1}\in A_{n+1}$ such that $L(b) \leq 0$. Fix $h\in \R^n$ such that $b_{n+1}=\sum_{i=1}^n h_i a_i$. For $t>0$ small enough, we may represent $z$ as:
\[z=\sum \lambda_i a_i -tb +tb=\sum_{i=1}^n (\lambda_i -t h_i) a_i + tb,\]
which is a valid representation for small $t>0$ , since $F_t:=\sum (\lambda_i-th_i)^q+t^q <1$  (using that $F(\lambda)<1$) and since $\lambda_i>0$ for all $i\in [n]$. And,
\begin{align*}
F_t=\sum (\lambda_i-th_i)^q+t^q &\geq F(\lambda) -qt\sum h_i \lambda_i^{q-1} +t^q \\
&=F(\lambda) -qt L(b) +t^q \geq F(\lambda) +t^q > F(\lambda),
\end{align*}
where we used that (for any $y\in \R$) $(1-ty)^q\geq 1-qty$ for $t>0$ small, and that $L(b)\leq 0$ by choice of $b\in A_{n+1}$. This contradicts maximality of $F(\lambda)$. Therefore we must have $r<n$: for that case as well, we shall get a contradiction.

So assume that $z=\sum_{i=1}^r \lambda_i a_i$, with $r\le n-1$, $a_1,...,a_r$ linearly independent, and $m(z)=F(\lambda)<1$. Denote $W=\text{span}(a_i)$ and define the linear form $L$ on $W$, by $L(\sum h_i a_i)=\sum h_i \lambda_i^{q-1}$. In particular, $L(z)=F(\lambda)>0$. Denote $H=\text{ker}(L)=\{w\in W : L(w)=0\}$. Hence $\dim(H)=r-1$, so $m:=\dim(H^{\perp})=n+1-r$. Denote $\pi$ the orthogonal projection onto $H^{\perp}$. Since $z\notin H$, $\pi z\neq 0$. By [Theorem 2.3, \cite{Barany82}], since $k-r \geq m$ and $0\in \text{conv}(\pi A_j)$ for $j\in [k]\setminus [r]$, we can pick $\alpha\geq 0$, $\mu_j\geq 0$ for $j=r+1,...,k$, and $b_j\in A_j$, such that $\alpha+\sum \mu_j=1$ and $\alpha \pi z+\sum \mu_j \pi b_j=0$. Set $M:=\sum_j \mu_j$, so $0<M\le 1$ (we cannot have $M=0$ because $\pi z\neq 0$). Denote $b=\sum \mu_j b_j$: since $b+\alpha z \in H$, we have $b\in -\alpha z+H \subset W$, so we can pick $h\in \R^r$ so to write $b=\sum h_i a_i$. And from $L(\alpha z+b)=0$ we also deduce that $L(b)\leq 0$.
Now, as in the case $r=n$, we use the newly picked $b_j$ to get other valid representations of $z$, namely :
\[z=\sum \lambda_i a_i -tb +tb=\sum_{i=1}^r (\lambda_i-th_i) a_i + \sum_{j=r+1}^k t\mu_j b_j \]
Denote $F_t := \sum_{i\leq r} (\lambda_i -th_i)^q +t^q \sum_{j>r} \mu_j^q$. Since $F(\lambda)=F_0 <1$, we have $F_t<1$ for $t>0$ small enough. Moreover
\begin{align*}
F_t= \sum_{i\leq r} (\lambda_i -th_i)^q +t^q \sum_{j>r} \mu_j^q &\ge F(\lambda) -qt \sum h_i \lambda_i^{q-1}+t^q \sum \mu_j^q \\
&=F(\lambda) -qt L(b) +t^q ||\mu||_q^q \\
& \ge F(\lambda) +t^q ||\mu||_q^q > F(\lambda)
\end{align*}
where we used that $L(b)\leq 0$, and that $||\mu||_q>0$ (since $||\mu||_1=M>0$). This contradicts the maximality of $F(\lambda)$ among valid representations.
We deduce that $m(z)=1$, and so $z\in R_p$, as desired.
\end{proof}


\begin{figure}[htbp]
    \centering
    \includegraphics[width=0.42\textwidth]{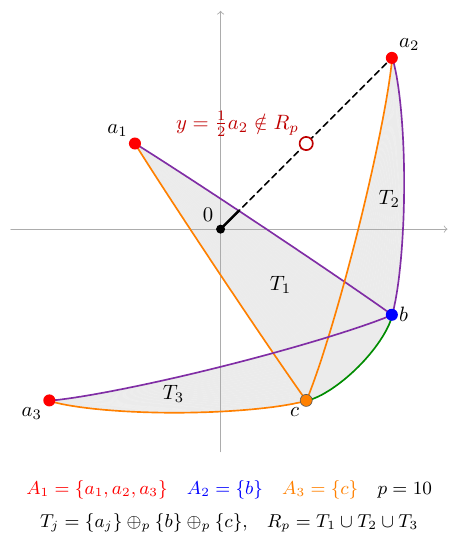}
   \caption{An example when $0\in \bigoplus_{i\in [3]}^{(\infty)} A_i$ but $R_p$ is not starshaped (at $0$).}
    \label{fig:contre-exemple}
\end{figure}

   Note that $n+1$ is optimal in the above statement. See remark \ref{rem:regularsimplex} below.

    \begin{remark}
    It is not enough to require that $0\in A_1 \oplus_{\infty} \cdots \oplus_{\infty} A_k$. Indeed, for $n \geq 2$, we give an example where one of the sets contains the origin and  $\bigoplus_{i \in [k]}^{(p)} A_i$ is not starshaped with respect to any point for $ p \in (1,\infty)$.
    For any $k\geq n + 1 $, we consider
    \[
        A_i = 
        \begin{cases}
                \{0,e_1 \} & i=1,
                \\
                \{e_i\} & 2\leq i\leq n
                \\
                \{e_n\} &i \geq n+1.
        \end{cases}
    \]
    By definition of $L_p$-sum, we have
    \begin{align}
        R_p: = \bigoplus_{i \in [k]}^{(p)} A_i 
        &= \left\{ \sum_{i=1}^k t_i^{1/q} a_i: a_i \in A_i, 0\leq t_i  \leq 1  \text{ and } \sum_{i=1}^k t_i =1\right\}.
    \end{align}
    First, notice that every $x \in R_p $ with $x_1 >0$, satisfies $\|x\|_q^q\geq 1$. Indeed, we may write $x \in R_p$ as  
    \[
        x = \sum_{i=1}^{n-1} t_i^{1/q}e_j + \left( \sum_{i = n}^{k} t_i^{1/q} \right) e_n, \qquad t_i \geq 0, \quad \sum_{i=1}^{k} t_i = 1.
    \]
    Then,
    \[
        \|x\|_q^q = \sum_{i=1}^{n-1} t_i + \left( \sum_{j = n}^{k} t_i^{1/q} \right)^q \geq \sum_{i=1}^{k} t_i = 1.
    \]
    Suppose, toward a contradiction, that $R_p$ is starshaped with respect to some $c \in R_p$. From the construction, all coordinates of $c$ are non-negative.
   If $c_1 > 0 $ then, since $0 \in R_p$, we have  $[0,c] \subset R_p$. For small enough $\lambda >0$, we have $\lambda c \in R_p$, but $\|\lambda c\|_q < 1$ and $\lambda c_1 >0$, which contradicts the above.
On the other hand, if $c_1 = 0$, then, since $e_1 \in R_p$, for any $t \in [0,1]$, we have $x = (1-t)e_1 + tc \in R_p.$ Thus,
        $\left\|x\right\|_q^q=(1-t)^q+t^q\|c\|_q^q.$
    Since $q >1$, for small enough $t>0$, we have $t^{q-1}\|c\|_q^q<1$, and so
    $$
        \left\|x\right\|_q^q<1-t+t\left(t^{q-1}\|c\|_q^q\right)<1.
    $$
    Since $x_1 = 1-t >0$, this again contradicts the above.
    Hence, $R_p $ is not starshaped with respect to any point, see the following picture in $\R^2$.
\end{remark}
\newcommand{\Pexp}{3}
\pgfmathsetmacro{\Qexp}{\Pexp/(\Pexp-1)}
\pgfmathsetmacro{\Lheight}{2^(1/\Pexp)}

\definecolor{Region}{HTML}{EAEAEA}

\newcommand{\MakeBand}{%
  \gdef\BandPath{}%
  \foreach \ii in {0,...,200}{%
    \pgfmathsetmacro{\ang}{90*(\ii/200)}%
    \pgfmathsetmacro{\xx}{max(cos(\ang),0)^(2/\Qexp)}%
    \pgfmathsetmacro{\yy}{\Lheight*max(sin(\ang),0)^(2/\Qexp)}%
    \xdef\BandPath{\BandPath (\xx,\yy)}%
  }%
  \foreach \ii in {0,...,200}{%
    \pgfmathsetmacro{\ang}{90-90*(\ii/200)}%
    \pgfmathsetmacro{\xx}{max(cos(\ang),0)^(2/\Qexp)}%
    \pgfmathsetmacro{\yy}{max(sin(\ang),0)^(2/\Qexp)}%
    \xdef\BandPath{\BandPath (\xx,\yy)}%
  }%
}

\MakeBand
\begin{center}
\begin{tikzpicture}[x=5.4cm,y=5.4cm,line cap=round,line join=round,
  every node/.style={font=\normalsize},scale = .4]
  \coordinate (O) at (0,0);
  \coordinate (Eone) at (1,0);
  \coordinate (Etwo) at (0,1);
  \coordinate (Top) at (0,\Lheight);
  \coordinate (Y) at (.5,0);

  \draw[gray!40,line width=1.3pt,-{Stealth[length=6pt]}]
    (-.16,0)--(1.20,0);
  \draw[gray!40,line width=1.3pt,-{Stealth[length=6pt]}]
    (0,-.16)--(0,1.4);

  \fill[Region] plot coordinates {\BandPath} -- cycle;

  \draw[black,line width=1.15pt] (O)--(Top);

  \draw[black,line width=.9pt]
    plot[domain=0:90,samples=201,variable=\t]
      ({max(cos(\t),0)^(2/\Qexp)},
       {max(sin(\t),0)^(2/\Qexp)});
  \draw[black,line width=.9pt]
    plot[domain=0:90,samples=201,variable=\t]
      ({max(cos(\t),0)^(2/\Qexp)},
       {\Lheight*max(sin(\t),0)^(2/\Qexp)});


  \fill[black] (O) circle[radius=2.6pt];
  \fill[black] (Eone) circle[radius=2.6pt];
  \node[below left=2pt] at (O) {$0$};
  \node[below right=2pt] at (Eone) {$e_1$};

  \path[fill=black] (Etwo) ++(90:2.8pt)
    arc[start angle=90,end angle=270,radius=2.8pt] -- cycle;
  \path[fill=black] (Etwo) ++(-90:2.8pt)
    arc[start angle=-90,end angle=90,radius=2.8pt] -- cycle;
  \node[left=4pt] at (Etwo) {$e_2$};

  \fill[black] (Top) circle[radius=1.7pt];
  \node[above left=0pt] at (0,1.2) {$2^{1/p}e_2$};

  \pgfmathsetmacro{\MidX}{(.5)^(1/\Qexp)}
  \pgfmathsetmacro{\MidY}{(1+\Lheight)*\MidX/2}

  \node[anchor=base] at (.50,-.50) {
    {$A_1=\{0,e_1\}$,}\quad
    {$A_2=\{e_2\}$,}\quad
    {$A_3=\{e_2\}$,}\quad $p=\Pexp$.};
\end{tikzpicture}
\end{center}
Here are two corollaries of Theorem \ref{thm:colorfulstarshapedness+prop11}.

\begin{corollary}
\label{coro:rainbowmonotonicity}
Let $k\geq n+1$ and let $A_1,...,A_k\subset \R^n$ be sets such that $0\in \cap_{i=1}^k \text{conv}(A_i)$. Then  $\left(\bigoplus_{i \in [k]}^{(p)} A_i\right)_p$
is non-increasing in $p>1$. In particular, 
$\bigoplus_{i \in [k]}^{(\infty)} A_i\subset \bigoplus_{i \in [k]}^{(p)} A_i$
for any $p>1$. If moreover $A_1,..., A_k$ are compact sets, then $\bigoplus_{i \in [k]}^{(p)} A_i$ converges to $\bigoplus_{i \in [k]}^{(\infty)} A_i$ as $p\to \infty$, for the Hausdorff metric.
\end{corollary}
\begin{proof}
Let $p'>p$ and let $z\in A_1\oplus_{p'} \cdots \oplus_{p'} A_k$. So $z=\sum \lambda_i a_i$, with $a_i\in A_i$, $\lambda_i\ge 0$ and $||\lambda||_{q'}=1$. So $\mu :=||\lambda ||^{-1}_q \ge 1$ (as $q'<q$), and $\mu z\in A_1\oplus_p \cdots \oplus_p A_k$. Since $z\in [0,\mu z]$ and $A_1\oplus_p \cdots \oplus_p A_k$ is starshaped, we deduce that $z\in A_1\oplus_p \cdots \oplus_p A_k$ as claimed.

 Now we turn to the moreover part, and assume all $A_i$ are compact. Let $R>0$ be such that $A_i\subset R B_2^n$ for all $i\in [k]$. If $\epsilon>0$ and $p>p(\epsilon, R)$, then $R_{\infty}\subset A_1\oplus_p \cdots \oplus_p A_k \subset R_{\infty}+\epsilon B_2^n$. Indeed, if $a_j\in A_i \subset R B_2^n$ ($j\le k$), and $z\in \{a_1\}\oplus_p \cdots\oplus_p \{a_k\}$, i.e. $z=\sum \lambda_i a_i$ with $||\lambda||_q=1$, then set $z'=\frac{z}{||\lambda||_1} \in R_{\infty}$, so $z\in ||\lambda||_1 R_{\infty}$, and $||\lambda||_1 \in [1, k^{1/p}]$. Therefore $z=z'+\left(||\lambda||_1-1\right) z'= z'+\epsilon u$ with $||u||_2=\epsilon^{-1} \left(||\lambda||_1-1\right) ||z'||_2 \le R \epsilon^{-1} (k^{1/p}-1) \le 1$ if $p \ge p(k,R, \epsilon)$ is large enough. 
\end{proof}

Taking $A_1=...=A_k=A$ in the above Corollary, one obtains the following statement. 
\begin{corollary}
\label{coro:starshapedloneA}
Let $A\subset \R^n$ be a set such that $0\in \text{conv}(A)$, and let $k\geq n+1$. Then $\text{conv}(A)\subset A^{\oplus_p k}$, and these sets are starshaped. More precisely, $(A^{\oplus_p k})_{p > 1}$ is non-increasing (for inclusion). If moreover $A$ is compact, then (for any $k\geq n+1$), $\lim_{p\to \infty} A^{\oplus_p k} =\text{conv}(A)$, in the Hausdorff metric.
\end{corollary}

\begin{remark}
\label{rem:regularsimplex}
We give an example which shows the sharpness of $n+1$ in Corollary \ref{coro:starshapedloneA} and in Theorem \ref{thm:colorfulstarshapedness+prop11}.
    Let $v_1, ... ,v_{n+1} \in \mathbb{S}^{n-1}$ be the vertices of a simplex $T$ such that $0\in\text{int}(T)$. Let $A=\{v_1, ... ,v_{n+1}\}$. Then $A^{\oplus_p k} \supset \text{conv}(A)$ for $k\geq n+1$, but $A^{\oplus_p n}$ doesn't contain the convex hull of $A$, as it doesn't contain $0$.
{\small {(indeed, $A^{\oplus_p n}=\cup_{j=1}^{n+1} B_j^{\oplus_p n}$, with $B_j=A\setminus \{v_j\}$. Say $B=B_{n+1}=\{v_1,...,v_n\}$. Then, since $0\in \text{int}(T)$, there exists an open half-space $H$, $H\supset B$ (and $0, v_{n+1}\notin H$). Hence $B^{\oplus_p n}\subset H^{\oplus_p n}=H$, implying $0\notin B^{\oplus_p n}$.)}} For $T$ a regular $n$-simplex, this argument gives $B(0,\frac{1}{n})\cap A^{\oplus_p n}=\emptyset$, while $A^{\oplus_p n} \cap \R_+ u\neq \emptyset$ for any $u\in \mathbb{S}^{n-1}$, so clearly $A^{\oplus_p n}$ isn't starshaped.
\end{remark}

\begin{remark}
The compactness assumption in the moreover part of the above corollary is needed. Indeed, consider the convex set $A=\{(x,y)\in \R^2; y\ge x^2\} \subset \R^2$. For any $k\ge 2$, and any $p\in (1,\infty)$, one has $A^{\oplus_p k}=k^{1/p}A$. And so $d_H(A^{\oplus_p k},A)=+\infty$ (since $(R,\frac{R^2}{k^{1/p}})\in k^{1/p} A$ and
{\small{$d_H\left(\left(R,\frac{R^2}{k^{1/p}}\right),A\right)\ge \frac{(1-k^{-1/p})R^2}{\sqrt{1+4R^2}}\xrightarrow[R\to\infty]{}+\infty$}}).
\end{remark}

\lineintersectionisconvex*

%
%
\begin{proof}
We can exclude the cases $A=\emptyset$ or $A=\{x_0\}$ (in which case $A^{\oplus_p k}=[x_0,k^{1/p}x_0]$), or the case $A\subset \R u_0$, for which the conclusion clearly holds.
If $A$ is not full-dimensional and $u\in \mathbb{S}^{n-1}\setminus \text{span}(A)$, then $\R u\cap A^{\oplus_p k} \subset \{0\}$ so there is nothing to show. Hence we can assume that $\text{span}(A)=\R^n$, and $n\ge 2$. Also, we can assume without loss of generality that $A$ is finite.
Either $0\in \text{conv}(A)$, then by Corollary \ref{coro:starshapedloneA}, $A^{\oplus_p k}$ is starshaped, in particular the intersection $I_u:=A^{\oplus_p k} \cap \R u$ is a segment (containing $0$). 

Or $0\notin \text{conv}(A)$. Since $A$ is finite, $\text{conv}(A)$ is compact, so by strict separation there exists $u_0\in\mathbb{S}^{n-1}$ and $\delta>0$ such that $\langle a,u_0\rangle \geq \delta >0$ for any  $a\in A$. In particular, denoting $H_0$ the open half-space $H_0=\{z\in \R^n : \langle z,u_0\rangle >0\}$, we see that for all $\ell \geq 1$, $A^{\oplus_p \ell}\subset H_0$.

We must show that for any $u\in \mathbb{S}^{n-1}$, $I_u:=\R_+ u\cap A^{\oplus_p k}$ is a segment. If $\langle u,u_0\rangle \leq 0$, then $\R_+ u\cap A^{\oplus_p k}=\emptyset$ (because $A^{\oplus_p k}\subset H_0$). Assume  $I_u=\R_+ u\cap A^{\oplus_p k}\neq \emptyset$, in particular $\langle u, u_0 \rangle >0$. Also assume $I_u$ has at least two points, else there is nothing to show.  Let $x\neq y$ be two arbitrary points within $I_u$. Hence (up to switching names between $x$ and $y$) there exists $0<r<s$ such that $x=ru$ and $y=su$. We must show that $[r,s]u\subset I_u$.

Fix $x=\sum_{i=1}^k \lambda_i a_i$ and $y=\sum_{j=1}^k \mu_j b_j$, two valid representations of $x$ and $y$, as points of $A^{\oplus_p k}$, i.e. with $a_i, b_j\in A$, $\lambda_i, \mu_j \geq 0$ and $||\lambda||_q=||\mu||_q=1$. Denote $(x_1,..., x_{2k})=(a_1,...,a_k,b_1,...,b_k)$ sequence of $2k$ vectors from $A$ (with perhaps some repeats). We define
\[P:=\left\{\nu\in \R^{2k} : \sum_{i=1}^{2k} \nu_i  x_i=u, \,\nu_j\geq 0 \, (\forall j\in [2k])\right\}\, .\]
If $\nu\in P$, then $\sum \nu_i \langle x_i, u_0\rangle=\langle u,u_0\rangle >0$, and since $\langle x_i,u_0\rangle \geq \delta$ for all $i$ (as $x_j\in A$), we deduce that $\sum_{i=1}^{2k} \nu_i \leq \delta^{-1} \langle u,u_0\rangle=:M >0$, implying in particular that $P\subset [0,M]^{2k}$ is compact. Thus $P$ is a polytope. The only inequalities used for defining $P$ are the $2k$ inequalities $\nu_i \geq 0$, so all faces of $P$ are sets of the form 
\[P_I:=\{\nu\in P : \text{supp}(\nu)\subset I\}=\{\nu \in P : i\in [2k]\setminus I \Rightarrow \nu_i=0\}=P\cap \left(\cap_{i\notin I} e_i^{\perp}\right) \]
for some $I$
(note that $P_I$ can be empty, depending on $I$). 

Denote $M=(x_1,...,x_{2k})$, which we see as an $n\times (2k)$-matrix. We claim:
\begin{claim}
If $P_I\neq \emptyset$, then $\text{dim}(P_I)=|I_0|-\text{rk}(M_{I_0})$ where $M_I$ is the $n\times I$-submatrix of $M$, i.e. $M_I=(x_i : i\in I)$, and where $I_0=\cup_{\nu \in P_I} \text{supp}(\nu)$ is the support of any $\nu \in \text{relint}(P_I)$.
\end{claim}
\begin{proof}[Proof of Claim]
Let us first argue in the case $I=I_0$.
So let $P_I\neq \emptyset$, let $\nu\in \text{relint}(P_I)$ and assume $\text{supp}(\nu)=I$. Then if $h\in \R^{2k}$ is such that $\nu \pm \epsilon h\in P$ for $\epsilon>0$ small enough, then necessarily $\text{supp}(h)\subset \text{supp}(\nu)\subset I$, and $\sum_{i\in I} h_i x_i=0$, i.e. $\tilde{h}\in \text{ker}(M_I)$, if $\tilde{h}=(h_i)_{i\in I}$ is obtained by deleting the (null) $I^c$-coordinates of $h$. Conversely, assume $h\in \R^{2k}$ with $\text{supp}(h)\subset I_0$, and with $\tilde{h}\in \text{ker}(M_I)$. Then because $\nu_i>0$ for all $i\in I$, we have $\nu\pm \epsilon h\in P_I$ when $\epsilon>0$ small enough. In other words, if we denote $\pi_I : \R^I \to \R^{2k}$ which gives back $h$ from $\tilde{h}$ (i.e. the embedding into $\R^{2k}$ given by $h_i=\tilde{h}_i$ for $i\in I$ and $h_i=0$ for $i\notin I$), then
\[\text{aff}(P_I)=\nu+\text{Span}(\pi_I(\text{ker}(M_I)))\]
and so $\dim(P_I)=\dim(\text{ker}(M_I))=|I|-\text{rk}(M_I)$, by the rank theorem.

In case $I_0\subsetneq I$, 
then $P_I=P_{I_0}$ is the same face of $P$ since they have the same relative interior. Thus the claim follows from the previous case.
\end{proof}
Denote $\mathcal{J}=\{I\subset [2k] : \exists \nu\in P, \text{supp}(\nu)=I\}$ (in other words, $I\in \mathcal{J}$ if for any $I'$ such that $P_{I'}=P_I$, one has $I\subset I'$).
By the above claim, if $I\in \mathcal{J}$ and $P_I$ is a $0$-dimensional face of $P$, then $|I|=\text{rk}(M_{I})\leq n$. Similarly,  if $I\in \mathcal{J}$ and $P_I$ is an edge of $P$, then  $|I|=\text{rk}(M_{I})+1\leq n+1 \leq k$. 
For $\ell \ge 0$, denote $\mathcal{J}_{\ell}:=\{I\in \mathcal{J} : \dim(P_I)=\ell\}$. Denote by $P^{(1)}$ the $1$-skeleton of $P$, i.e. the union of its edges (and vertices):
\[P^{(1)}=\bigcup_{I\in \mathcal{J}_0\cup \mathcal{J}_1} P_I=\bigcup_{I\in\mathcal{J}_1} P_I\]
(where the second equality is due to $\dim(P)\ge 1$).
All faces of $P$ are of type $P_{I_0}$ for some $I_0\in \mathcal{J}$. We deduce that the $1$-skeleton of $P$ is contained inside
\[P_{\leq k} :=\{\nu\in P : |\text{supp}(\nu)|\leq k\}=\bigcup_{I\subset [2k], |I| \leq k} P_{I}=\bigcup_{I_0\in \mathcal{J}, |I_0| \leq k} P_{I_0} \hspace{2mm}. \]
We claim that $P_{\leq k}$ is path-connected.
Indeed, if $\nu^{(0)},\nu^{(1)}$ lie within $P_{\leq k}$, let $I_0, I_1$ denote their respective support, thus $\nu^{(0)}\in P_{I_0}$ and $\nu^{(0)}\in P_{I_1}$. For each of these two faces of $P$, pick a vertex $\mu^{(j)}\in \text{Vert}(P_{I_j})$, $j=0,1$. The $\mu^{(j)}$ are vertices of $P$ (possibly the same one), hence there is a path on $P^{(1)}$ (a path made of edges $(e_i)_{1\le i\le m}$ with $e_i$, $e_{i+1}$ sharing a vertex) linking $\mu^{(0)}$ to $\mu^{(1)}$. And obviously (for $j=0,1$) the  segments $[\mu^{(j)},\nu^{(j)}]\subset P_{I_j} \subset P_{\leq k}$. Concatening these two segments with the path on $P^{(1)}\subset P_{\leq k}$, yields a continuous path from $\nu^{(0)}$ to $\nu^{(1)}$, within $P_{\leq k}$, showing path-connectedness of the latter.

We can now conclude the proof that $[x,y]\subset I_u=\R_+ u\cap A^{\oplus_p k}$. Since $x=ru=\sum_{i=1}^k \lambda_i a_i$, we see that $\nu^{(0)}:=(\frac{\lambda}{r},0)=(\frac{\lambda_1}{r},...,\frac{\lambda_k}{r},0,...,0) \in P_{\leq k}$. Similarly, since $y=su=\sum_{j=1}^k \mu_j b_j$, we see that $\nu^{(1)}:=(0,\frac{\mu}{s})\in P_{\leq k}$. Since $P_{\leq k}$ is path-connected, there exists a continuous $\gamma : [0,1] \to P_{\leq k}$, with $\gamma(0)=\nu^{(0)}$ and $\gamma(1)=\nu^{(1)}$. Note that $||\gamma(t)||_q >0$ for all $t\in [0,1]$, since $\nu=0_{\R^{2k}}\notin P$. Thus the function $f(t):=\frac{1}{||\gamma(t)||_q}$ is continuous on $[0,1]$. Observe that $f(0)=r$ and that $f(1)=s$ (since $||\lambda||_q=||\mu||_q=1$). By continuity, $f([0,1])$ is a segment, so $[r,s]\subset f([0,1])$. Hence for any $r'\in [r,s]$, there exists $t\in [0,1]$ such that $f(t)=r'$. Denote $\nu=\gamma(t)$, then $||r' \nu||_q=1$ and $\sum_{i\in [2k]} (r'\nu_i) x_i=r' u$ by definition of $P$. Moreover $|\text{supp}(r' \nu)|=|\text{supp}(\nu)|\leq k$ since $\gamma([0,1])\subset P_{\le k}$, and so $r' u\in A^{\oplus_p k} \cap \R u$.
We conclude that $[x,y]=[r,s]u\subset I_u$, as desired.
\end{proof}

 For $B\subset \R^n$ a measurable set, let us use the notation $B^*=\cup_{b\in B} [0,b]$.
 
\begin{lemma}
\label{lem:geomrad}
Suppose that $B\cap\R u$ is an interval for every \(u\in\mathbb S^{n-1}\). Then either $B=B^*$, or $B^*\setminus B$ is starshaped.
\end{lemma}

\begin{proof}
If \(0\in B\), the assumption implies that \([0,y]\subset B\) for every \(y\in B\), and so $B=B^*$. Now assume $0\notin B$, in particular $0\in B^*\setminus B$.
Let $x\in B^*\setminus B$. We can write \(x=ry\) for some \(y\in B\) and $r\in [0,1)$. Note that our assumption on $B$, and that $y\in B$, $x\notin B$, implies that $[0,x] \cap B=\emptyset$. So clearly $[0,x] \subset B^*\setminus B$.
\end{proof}

\begin{corollary} \label{cor:subtraction-starshaped}
    Let $1<p\le\infty$ and let $A \subset \mathbb{R}^n$ be a measurable set. Then, for every $k$, one has $(A^*)^{\oplus_p k}=(A^{\oplus_p k})^*$. For $k \geq n+1$, the set $S_k:=(A^*)^{\oplus_p k} \setminus A^{\oplus_p k}$ is either empty or starshaped. 
    Moreover, $S_k = S_{n+1} \subset S_n$ for all $k \geq n+1$, and so $A^{\oplus_p k}=\left(A^{\oplus_p k}\right)^* \setminus S_{n+1}$.
\end{corollary}

\begin{proof}
We first show that $(A^*)^{\oplus_p k}=(A^{\oplus_p k})^*.$ Observe that $(A^*)^{\oplus_p k}$ is starshaped (for any $k\ge 1$). 
As $A^{\oplus_p k}\subset (A^*)^{\oplus_p k}$, we have $(A^{\oplus_p k})^*\subset ((A^*)^{\oplus_p k})^*=(A^*)^{\oplus_p k}$. 
Conversely, if $x\in (A^*)^{\oplus_p k}$, then we can write
$x=\sum_{i=1}^k\lambda_i \rho_i a_i$,
 for some 
$a_i\in A$, $\rho_i\in[0,1]$,
$\lambda_i\geq0$, $\|\lambda\|_q=1.$
Set $\rho=\|(\lambda_i \rho_i)_{i}\|_q\leq1$.
If $\rho=0$, then $x=0$. Otherwise,
$x=\rho\cdot \sum_{i=1}^k\frac{\lambda_i \rho_i}{\rho}a_i=: \rho\cdot y$, with $y\in A^{\oplus_p k}$. 
This shows the reverse inclusion.

Now assume $k\ge n+1$. By Theorem~\ref{prop:lineintersectionisconvex},
$A^{\oplus_p k}\cap\mathbb Ru$ is an interval for every
$u\in\mathbb S^{n-1}$.
Applying Lemma \ref{lem:geomrad}, we deduce that $(A^{\oplus_p k})^*\setminus A^{\oplus_p k}$ is either empty, or starshaped.


Next, we show that $S_{k+1}=S_k\subset S_n$ if $k\ge n+1$. This is clear if $0\in \text{conv}(A)$, since $S_k=\emptyset$ for $k\ge n+1$, in that case. So we assume that $0\notin \text{conv}(A)$.
First, let us show that $S_{k+1}\subset S_k$, for $k\ge n$.   Let $x\in (A^*)^{\oplus_p (k+1)} \setminus A^{\oplus_p (k+1)}$. Clearly $x\notin A^{\oplus_p k}$, since $A^{\oplus_p k} \subset A^{\oplus_p (k+1)}$. Let us show that $x\in (A^*)^{\oplus_p k}$. We can assume $x\neq 0$. Since $x\in \text{cone}(A)$, by Carath\'eodory's theorem we can write $x=\sum_{i=1}^m \lambda_i a_i$, for some $1\le m\le \min\{k,n\}\le k$, some $\lambda_i\ge 0$ and some $a_i\in A$. Call $x'=\frac{x}{||\lambda||_q}$. Thus $x'\in A^{\oplus_p m} \subset  A^{\oplus_p k}$. If we had $||\lambda||_q \ge 1$, then we would have $x'\in [0,x]$, and since $x\in (A^*)^{\oplus_p (k+1)} \setminus A^{\oplus_p (k+1)}$ and the latter is starshaped, this would imply that $x'$ lies in this set, in particular $x'\notin A^{\oplus_p (k+1)}$ (a contradiction). Hence $||\lambda||_q<1$ and $x\in [0,x']$. Since $x'\in A^{\oplus k}$, we deduce that $x\in (A^{\oplus_p k})^*$, and hence $S_{k+1}\subset S_k$ for $k\ge n$.  

 If $p=\infty$, then $S_k=\text{conv}(A^*)\setminus \text{conv}(A)$ for all $k\ge n+1$, by Carath\'eodory theorem. Similarly, let us argue that $S_k=S_{n+1}$ for all $k\ge n+1$, when $p\in (1,\infty)$, which will close the proof of Corollary \ref{cor:subtraction-starshaped}. Thus we only need to check that $S_{n+1}\subset S_{k}$, if $k\ge n+2$. Since $(A^*)^{\oplus_p (n+1)}\subset (A^*)^{\oplus_p k}$, it suffices to show that if $x\in S_{n+1}$, then $x\notin A^{\oplus_p k}$. 
So assume for contradiction there was some $x\in S_{n+1} \cap A^{\oplus_p k}$, and recall we are assuming $0\notin \text{conv}(A)$.
 We write $x=\sum_{i=1}^k \lambda_i a_i$, for some $a_i\in A$, $\lambda_i\ge 0$, with $||\lambda||_q=1$.
 By Lemma \ref{lem:coefficient-reduction-q}, we can find $I\subset [k]$ of size $s\leq n$, and $c\in \R_+^I$, such that $x=\sum_{i\in I} c_i a_i$, $c_i\ge 0$, and $||c||_q\ge ||\lambda||_q=1$. Then $x'=||c||_q^{-1} x \in A^{\oplus_p n}\subset A^{\oplus_p (n+1)}$ ;
and $x'=||c||_q^{-1} x\in [0,x]$. Since $x\in S_{n+1}$ which is starshaped, we deduce that $x'\in S_{n+1}$ as well (contradiction).
Consequently, $S_{n+1}\subset S_k$ for every $k\ge n+1$. Together with the reverse inclusion
established above, this yields that 
$S_n\supset S_{n+1}=S_k$, for all $k\ge n+1$.
 \end{proof}
\begin{remark}
     For $p =\infty$, the inclusion $S_{n+1}\subset S_n$ says $\text{conv}(A^*)\setminus \text{conv}(A) \subset (A^*)^{\oplus_{\infty} n} \setminus A^{\oplus_{\infty} n}$. We note that this inclusion can be strict : if $A=\text{Vert}(T)$ with $T$ an $n$-simplex, $0\in \text{int}(T)$, then $\text{conv}(A^*)\setminus \text{conv}(A)=\emptyset$ while $(A^*)^{\oplus_{\infty} n} \setminus A^{\oplus_{\infty} n}=\text{int}(T)$. For $p\in (1,\infty)$, for the same $A$, observe that $A^{\oplus_p (n+1)}$ is starshaped (by Corollary \ref{coro:starshapedloneA}), while $\text{int}(T)\subset (A^*)^{\oplus_p n} \setminus A^{\oplus_p n}$, so $\emptyset=S_{n+1}\subset S_n$ is strict too.  More precisely, for that example, we have $(A^*)^{\oplus_p n} \setminus A^{\oplus_p n}=\text{int}(T_p)$, where $T_p=\text{conv}\left(\cup_{i=1}^{n+1} \bigoplus^{(p)}_{j\in [n+1]\setminus \{i\}} \{a_j\}\right).$ 
\end{remark}

\section{Hausdorff metric}

In this section, we prove Theorem \ref{thm:Firey_convergence} and we establish sharp upper bounds for the Hausdorff distance of the $L_p$-sum of two sets. In fact, Theorem \ref{thm:Firey_convergence} is a consequence of Theorems \ref{thm: hausd-convergence-2}, \ref{theorem:lp_hausdorff_dimension_2_max_bound}, \ref{theorem:hausdorff_p_subadditive}, and \ref{theorem:strong_fractional_subadditivity}, which are proved below. Each of these theorems approaches the convexification phenomenon from a different perspective. We begin with a proof that relies heavily on the geometric properties of the $L_p$-sum that we derived in Section \ref{section:goemetry_of_lp_sums}.

 \begin{theorem}
 \label{thm: hausd-convergence-2}  Let $A \subset \mathbb{R}^n$ be a nonempty compact set, and let $1<p<\infty$. Then, for every $k \geq n+1$,
    \begin{equation} \label{eq:d-zero-convergence}
        d_0\left(A^{\oplus_p k}\right) \leq \max \left\{\min _{a \in A}\|a\|_2,\left(\max _{a \in A}\|a\|_2\right)\left[k^{1 / p}-(k-n+1)^{1 / p}\right]\right\} .
    \end{equation}
    In particular, for sufficiently large $k$,
    \[
        d_H\left(\frac{1}{k^{1/p}}A^{\oplus_p k}, \conv(0,A)\right) \leq  \frac{d_0(A)}{k^{1/p}}.
    \]
    Moreover, the above bound is sharp.
     \end{theorem}
     \begin{proof}
        Let $ k \geq n+1$.
        Using Lemma \ref{lem:containing-hyperplanebis}, we get
        \begin{equation} \label{eq:shell-standard}
            \left\{(x_1,\ldots,x_n)\in \R_+^n: \| x \|_1 \leq (k-n+1)^{1/p}, \|x\|_q \geq 1 \right\} \subset  \{e_1,\ldots,e_n\}^{\oplus_p k}.
        \end{equation}
        Using Lemma \ref{lem:coefficient-reduction}, every point $x \in \conv (0,A)$ can be written as
        \[
            x=\sum_{i=1}^n c_i a_i, \quad a_i \in A, \quad c_i \geq 0, \quad \sum_{i=1}^n c_i \leq 1.
        \]
        Define the linear map $T: \mathbb{R}^{n} \longrightarrow \mathbb{R}^n$ such that $T(e_j) = a_j$ for $1\leq j \leq n$. Set $u = k^{1/p}c$.
        Then, 
        $$
            T(u) = k^{1/p} \sum_{i=1}^n c_iT(e_i) = k^{1/p}x, \quad \|u\|_1 = k^{1/p}\|c\|_1 \leq k^{1/p}.
        $$
        Since
        $$
            T(u) \in k^{1/p}\conv (0,A) = k^{1/p} \conv \left( 0, \frac{A^{\oplus_p k}}{k^{1/p}} \right) = \conv \left(0,A^{\oplus_p k}\right),
        $$
        the desired inequality \eqref{eq:d-zero-convergence} follows from the following claim
        \[
            d\left(T(u), A^{\oplus_p k}\right) \leq \max \left\{\min _{a \in A}\|a\|_2,\left(\max _{a \in A}\|a\|_2\right)\left[k^{1 / p}-(k-n+1)^{1 / p}\right]\right\}.
        \]
        Note that if $v \in\left\{e_1, \ldots, e_n\right\}^{\oplus_p k}$, then
        \begin{equation} \label{eq:Tv-in-ei}
            T(v) \in T\left(\left\{e_1, \ldots, e_n\right\}^{\oplus_p k}\right)=\left\{a_1, \ldots, a_n\right\}^{\oplus_p k} \subset A^{\oplus_p k} .    
        \end{equation}
            \begin{enumerate}[wide =0pt, labelwidth= 1.1cm,leftmargin=1.55cm, ]
    \item[Case 1:] Suppose $\|u\|_q < 1$. Let $a_0$ be a point in $A$ such that $\|a_0\|_2 = \min_{a\in A} \|a\|_2$ and set $ t = \left(1-\|u\|_q^q\right)^{1 / q}$. Since $ k \geq n+1$ and by the definition of $L_p$-sum,
    \[
        T(u) +ta_0 = \sum_{i=1}^n u_i a_i  + ta_0 \in A^{\oplus_p (n+1)} \subset A^{\oplus_p k}.
    \]
    Thus,
    \[
        d (T(u),  A^{\oplus_p k}) \leq t \|a_0\|_2 \leq \min _{a \in A}\|a\|_2.
    \]
    \item[Case 2:] If $\|u\|_q \geq  1$ and $\|u\|_1 \leq (k-n+1)^{1 / p}$, then it follows from \eqref{eq:shell-standard} that $ u \in \left\{e_1, \ldots, e_n\right\}^{\oplus_p k}$ and so from \eqref{eq:Tv-in-ei}, $T(u) \in A^{\oplus_p k}$ and so $d(T(u),  A^{\oplus_p k}) = 0$.
    \item[Case 3:] If $\|u\|_q \geq  1$ and $\|u\|_1 > (k-n+1)^{1 / p}$, we take
    \[
        v = \max \left\{ \frac{(k-n+1)^{1 / p}}{\|u\|_1}, \frac{1}{\|u\|_q} \right\} u.
    \]
    Then, $\|v\|_q \geq 1$. If the maximum is attained by the first term, then $\|v\|_1=(k-n+1)^{1/p}$. On the other hand, if it is attained by the second term, then $\|v\|_q=1$.
    Thus, by \eqref{eq:shell-standard} in the first case and Lemma \ref{lem:sum-of-ei} in the second case, we obtain in either case that $v \in \left\{e_1, \ldots, e_n\right\}^{\oplus_p k}$. It follows from \eqref{eq:Tv-in-ei} that $T(v) \in A^{\oplus_p k}$. Since $\|u\|_1 \leq k^{1/p}$, 
    \begin{align}
        \|u-v\|_1 &= \left( 1- \max \left\{ \frac{(k-n+1)^{1 / p}}{\|u\|_1}, \frac{1}{\|u\|_q} \right\} \right)\|u\|_1 \\
        & \leq \|u\|_1 - (k-n+1)^{1 / p}\leq k^{1/p} - (k-n+1)^{1 / p}.
    \end{align}
    By the triangle inequality, we have
    \begin{align}
        d(T(u), A^{\oplus_p k}) &\leq  \| T(u) - T(v) \|_2 =\left\|\sum_{i=1}^n\left(u_i-v_i\right) a_i\right\|_2
        \\
        &\leq \sum_{i=1}^n\left|u_i-v_i\right|\left\|a_i\right\|_2 \leq \left(\max _{a \in A}\|a\|_2\right) \|u - v\|_1 
        \\
        &\leq \left(\max _{a \in A}\|a\|_2\right)\left[k^{1 / p}-(k-n+1)^{1 / p}\right].
    \end{align}
    \end{enumerate} 
    Combining all cases proves the claim. Taking the sup over $x \in \conv (0,A)$ yields \eqref{eq:d-zero-convergence}.

    Finally, if $d_0 (A) =0$, then $A $ is convex containing the origin, and so $A^{\oplus_p k}=k^{1 / p} A$. Thus, $$d_0\left(\frac{A^{\oplus_p k}}{k^{1/p}}\right) = 0.$$ If $d_0 (A) >0$, then for sufficiently large $k$ both terms in the maximum are at most $d_0 (A)$. Indeed, this follows from the fact that 
    \[
     \min_{a \in A} \|a\|_2 = d(0,A) \leq d_0(A),
    \]
    and for $p >1$,
    \[
        k^{1/p} - (k-n+1)^{1/p} \to 0, \qquad k \to \infty. 
    \]
Finally, we consider $A = \{a\}$ with $a \neq 0$. Then, $\frac{1}{k^{1/p}} A^{\oplus_p k} = [k^{-1/p}a,a]$ and $\conv (0,A) = [0,a]$. Thus,
    \[
        d_H\left(k^{-1 / p} A^{\oplus_p k}, \conv (0,A)\right)=\frac{\|a\|_2}{k^{1 / p}} = \frac{d_0(A)}{k^{1/p}},
    \]
    verifying the sharpness.
     \end{proof}

Recall the following sharp upper bound for the sum of two planar sets that we stated in the introduction.

\HausdorffDimTwoBound*

We state the first key lemma in the proof of Theorem \ref{theorem:lp_hausdorff_dimension_2_max_bound}, which was originally proved by the third-named author in the papers \cite{M-2025, M-2026}.

\begin{lemma}\label{lemma:parallelogram_lemma}
    Let $\tilde{A},\tilde{B}\subset\mathbb{R}^2$ be nonempty compact sets. Let $z\in \textup{conv}(\tilde{A}+\tilde{B})$ such that 
    \begin{align*}
        z\notin \tilde{A}+\textup{conv}(\tilde{B}),\qquad z\notin \textup{conv}(\tilde{A})+\tilde{B}.
    \end{align*}
    Then there exist distinct $a_0,a_1\in \tilde{A}$ and distinct $b_0,b_1\in \tilde{B}$ such that 
    \begin{align*}
        z\in [a_0,a_1]+[b_0,b_1],
    \end{align*}
    and that 
    \begin{align*}
        |a_0-a_1|\leq 2d(\tilde{A}),\qquad |b_0-b_1|\leq 2d(\tilde{B}).
    \end{align*}
    Moreover, the intervals $[a_0,a_1]$ and $[b_0,b_1]$ are nonparallel, so that the set $[a_0,a_1]+[b_0,b_1]$ is a parallelogram whose vertices are contained in $\tilde{A}+\tilde{B}$.
\end{lemma}

In conjunction with Lemma \ref{lemma:parallelogram_lemma}, we will use  Lemma \ref{lemma:firey_parallelogram_estimate} below, which we shall deduce from the following lemma.

\begin{lemma}\label{lemma:generalised_parallelograms}
    For $q\in [1,\infty)$, let $A,B\subset \R^n$ be two compact convex sets such that $1\leq \dim(A),\dim(B)\leq n-1$. Define $P_t=(1-t)^{1/q} A+t^{1/q} B$. Then if $z\in \cup_{t\in [0,1]} P_t$, then $z\in \partial P_{t_1}$ for some $t_1\in [0,1]$.
\end{lemma}
\begin{proof}
If $z\in A\cup B=\partial P_0 \cup \partial P_1$, there is nothing to show : hence we may assume $z\notin A\cup  B$. 


For $u\in \mathbb{S}^{n-1}$, denote $h_t(u)=(1-t)^{1/q} h_{A}(u)+t^{1/q} h_B(u)$. Therefore
\[P_t=(1-t)^{1/q} A+t^{1/q} B=\{z\in \R^n  :\forall u\in \mathbb{S}^{n-1}, \langle z,u\rangle \leq h_t(u)\}.\]
Since $z\notin B$, there exists $u_0\in \mathbb{S}^{n-1}$ such that $\langle z,u_0\rangle > h_B(u_0)$. By continuity, it implies that $\langle z,u_0\rangle > h_t(u_0)$ (thus $z\notin P_t$) for $t$ close enough to $1$.
 This shows that $t_1<1$, where
\[t_1=\sup \{t\geq 0 : z\in P_t\}.\]
If we had $z\in \text{int}(P_{t_1})$, then by compactness we would have $\min_{u\in\mathbb{S}^{n-1}} (h_{t_1}(u)-\langle z,u\rangle) >0$, which would imply that $z\in \text{int}(P_t)$ for $t>t_1$ close enough to $t_1$, contradicting the definition of $t_1$. Therefore $z\in \partial P_{t_1}$.
\end{proof}

As a particular case of Lemma \ref{lemma:generalised_parallelograms}, we find:
\begin{lemma}\label{lemma:firey_parallelogram_estimate}
    Let $I=[a_0,a_1]$ and $J=[b_0,b_1]$ be nonparallel segments in $\mathbb{R}^2$ such that 
    \begin{align*}
        |a_0-a_1|\leq 2r,\qquad |b_0-b_1|\leq 2r,
    \end{align*}
    for some $r>0$. For $t\in [0,1]$, define the parallelogram
    \begin{align*}
        P_t:=(1-t)^{1/q}I+t^{1/q}J,
    \end{align*}
    and define the vertices of $P_t$ by
    \begin{align*}
        \gamma_{ij}(t):=(1-t)^{1/q}a_i+t^{1/q}b_j,\qquad i,j\in\{0,1\}.
    \end{align*}
    If $z\in \cup_{t\in [0,1]} P_t$, then there exists $t_1\in [0,1]$ such that $z\in \partial P_{t_1}$. In particular,
    \begin{align*}
        d(z,\{\gamma_{ij}(t):0\leq t\leq 1,i,j\in\{0,1\}\})\leq r.
    \end{align*}
\end{lemma}
\begin{proof}
Apply Lemma \ref{lemma:generalised_parallelograms} with $A=I$, $B=J$, $n=2$. Since $I$ and $J$ are not parallel, we indeed have $\dim(I+J)=2$. If $z\in I\cup J$, the conclusion obviously holds. Else, we found $t_1\in (0,1)$ such that $z\in \partial P_{t_1}$, we deduce that $z$ lies on a translate of $\alpha I$, or on a translate of $\beta J$, with $\alpha=(1-t_1)^{1/q}\in (0,1)$ and $\beta=t_1^{1/q}\in (0,1)$. Hence there exists $c_{t_1}=\gamma_{ij}(t_1)$ ($i,j\in \{0,1\}$) a vertex of $P_{t_1}$, such that $|z-c_{t_1}|\leq \max\{\alpha r, \beta r\} \leq r$, as claimed.
\end{proof}



\begin{proof}[Proof of Theorem \ref{theorem:lp_hausdorff_dimension_2_max_bound}]
    Let $z\in \textup{conv}(0,A\oplus_p B)$. Recall that, from Lemma~\ref{lemma:conv0p} one has
    \begin{align*}
        \textup{conv}(0,A\oplus_p B)=\textup{conv}(0,A)\oplus_p\textup{conv}(0,B),
    \end{align*}
    so $z\in \textup{conv}(0,A)\oplus _p \textup{conv}(0,B)$. If $z\in (\{0\}\cup A)\oplus_p \textup{conv}(0,B)$, then we write $z=t_z^{1/q}a+(1-t_z)^{1/q}b$, for $a\in \{0\}\cup A$, $b\in \textup{conv}(0,B)$, and $0\leq t_z\leq 1$. If $a\in A$, then we have
    \begin{align*}
        d(z,A\oplus_p B)&\leq d(z,t_z^{1/q}a+(1-t_z)^{1/q}B)
        =(1-t_z)^{1/q}d(b,B)
        \leq \max\{d_0(A),d_0(B)\}.
    \end{align*}
    If $a=0$, then $z=(1-t_z)^{1/q}b$, and so we have $z\in \textup{conv}(0,B)$. Then
    \begin{align*}
        d(z,A\oplus_p B)\leq d(z,B)\leq \max\{d_0(A),d_0(B)\}.
    \end{align*}
    Alternatively, if $z\in \textup{conv}(0,A)\oplus _p (\{0\}\cup B)$, then we use the same technique to show that $d(z,A\oplus _p B)$ is bounded above by $\max\{d_0(A),d_0(B)\}$.

    Next, we consider the case where 
    \begin{align}\label{eq:non_containment_lp_condition}
        z\notin \textup{conv}(0,A)\oplus_p (\{0\}\cup B),\qquad z\notin (\{0\}\cup A)\oplus_p \textup{conv}(0,B).
    \end{align}
    Since $z\in \textup{conv}(0,A)\oplus_p\textup{conv}(0,B)$, there exists $0< t_z< 1$ such that $z\in t_z^{1/q}\textup{conv}(0,A)+(1-t_z)^{1/q}\textup{conv}(0,B)$. Moreover, by condition \eqref{eq:non_containment_lp_condition}, we have
    \begin{align*}
        z\notin t_z^{1/q}(\{0\}\cup A)+(1-t_z)^{1/q}\textup{conv}(0,B),\qquad z\notin t_z^{1/q}\textup{conv}(0,A)+(1-t_z)^{1/q}(\{0\}\cup B).
    \end{align*}
    It follows from Lemma \ref{lemma:parallelogram_lemma} with $\tilde{A}:=t_z^{1/q}(\{0\}\cup A)$ and $\tilde{B}:=(1-t_z)^{1/q}(\{0\}\cup B)$ that there exist $a_0,a_1\in \{0\}\cup A$ and $b_0,b_1\in \{0\}\cup B$ such that the segments $[a_0,a_1]$ and $[b_0,b_1]$ are nonparallel, that 
    \begin{align*}
        z\in t_z^{1/q}[a_0,a_1]+(1-t_z)^{1/q}[b_0,b_1],
    \end{align*}
    and that 
    \begin{align*}
        |a_0-a_1|\leq 2d(\{0\}\cup A),\qquad |b_0-b_1|\leq 2d(\{0\}\cup B).
    \end{align*}
    But notice that 
    \begin{align*}
        d(\{0\}\cup A)\leq d_0(A),\qquad d(\{0\}\cup B)\leq d_0(B),
    \end{align*}
    so the distances $|a_0-a_1|$ and $|b_0-b_1|$ are bounded above by $2\max\{d_0(A),d_0(B)\}$.

    If none of the points $a_0,a_1,b_0,b_1$ is the origin, then apply Lemma \ref{lemma:firey_parallelogram_estimate} to $I:=[a_0,a_1]$ and $J:=[b_0,b_1]$ with $r:=\max\{d_0(A),d_0(B)\}$ to get
    \begin{align*}
        d(z,A\oplus_p B)\leq d(z,\{\gamma_{ij}(t):0\leq t\leq 1,i,j\in\{0,1\}\})\leq \max\{d_0(A),d_0(B)\}.
    \end{align*}
    It remains to consider the case where an endpoint is zero. Suppose that $[a_0,a_1]=[0,a]$ with $a\in A$. Then there exists $0\leq \gamma\leq 1$ and $b\in \textup{conv}(0,B)$ such that 
    \begin{align*}
        z=t_z^{1/q}\gamma^{1/q} a+(1-t_z)^{1/q}b.
    \end{align*}
    Set $\tilde{t}_z:=t_z\gamma$. If $\tilde{t}_z<1$, then set
    \begin{align*}
        \tilde{b}:=\left(\frac{1-t_z}{1-\tilde{t}_z}\right)^{1/q}b.
    \end{align*}
    In fact, $\tilde{b}\in \textup{conv}(0,B)$. Then 
    \begin{align*}
        z=\tilde{t}_z^{1/q}a+(1-\tilde{t}_z)^{1/q}\tilde{b}
    \end{align*}
    We have already considered such a point $z$ at the beginning of this proof, and we found that $d(z,A\oplus _p B)$ is bounded above by $\max\{d_0(A),d_0(B)\}$. If $\tilde{t}_z=1$, then $z=a\in A\oplus _p B$ and $d(z,A\oplus_p B)=0$. The same technique works if any of the other endpoints $a_1,b_0,b_1$ is $0$.

    We showed above that if $z\in \textup{conv}(0,A\oplus_p B)$, then
    \begin{align*}
        d(z,A\oplus_p B)\leq \max\{d_0(A),d_0(B)\},
    \end{align*}
    and by taking the supremum over all such $z$, we have the desired bound.

    Next, we show that when $p=\infty$, the maximum bound holds for the measure $d$. Let $x\in \textup{conv}(A\oplus_{\infty} B)$. Then there exists $t\in [0,1]$ such that $x\in \textup{conv}(tA)+\textup{conv}((1-t)B)$. Similar to above, if $x\in tA+\textup{conv}((1-t)B)$, or if $x\in \textup{conv}(tA)+(1-t)B$, then
    \begin{align*}
        d(x,A\oplus_{\infty}B)\leq \max(d(A),d(B)).
    \end{align*}
    Otherwise, by Lemma \ref{lemma:parallelogram_lemma} there exist distinct $a_1,a_2\in A$ with $\|a_1-a_2\|_2\leq 2d(A)$ and distinct $b_1,b_2\in B$ with $\|b_1-b_2\|_2\leq 2d(B)$ such that
    \begin{align*}
        x\in t[a_1,a_2]+(1-t)[b_1,b_2]\subset [a_1,a_2]\oplus_{\infty}[b_1,b_2].
    \end{align*}
    Note the formula
    \begin{align*}
        [a_1,a_2]\oplus_{\infty}[b_1,b_2]=\bigcup_{j=1}^{2}[\textup{conv}([a_1,a_2],b_j)\cup \textup{conv}([b_1,b_2],a_j)],
    \end{align*}
    whose verification is an exercise in planar geometry. Suppose that $x\in \textup{conv}([a_1,a_2],b_j)$ for some $j$. By definition, we have $[b_j,a_1]\cup [b_j,a_2]\subset A\oplus_{\infty}B$. It follows immediately that $b_j,a_1,a_2$ are not collinear, since otherwise $x\in A\oplus_{\infty}B$. Then $x$ is trapped in the triangle $\textup{conv}([a_1,a_2],b_j)$, and because the sides $[b_j,a_1]$, $[b_j,a_2]$ are contained in $A\oplus_{\infty}B$, the distance from $x$ to the set $A\oplus_{\infty}B$ is at most half the length of the interval $[a_1,a_2]$, and the length of this interval is at most $2d(A)$. That is, $d(x,A\oplus_{\infty}B)\leq d(A)$. If $x\in \textup{conv}([b_1,b_2],a_j)$ for some $j$, then the same reasoning gives $d(x,A\oplus_{\infty}B)\leq d(B)$. In any case, $d(x,A\oplus_{\infty}B)$ is bounded above by $\max(d(A),d(B))$.

     To see that equality is possible, even if both sets are nonconvex, take the sets
    \begin{align*}
        A:=[0,2\alpha e_1]\cup [0,2\alpha e_2],\qquad B:=[0,-2\beta e_1]\cup [0,-2\beta e_2].
    \end{align*}
    Then $d_0(A)=\alpha$ and $d_0(B)=\beta$, and we easily deduce that $d_0(A\oplus_p B)=\max\{\alpha,\beta\}$.
\end{proof}

Next, we prove a sharp two-set bound for the case $n\geq 3$.

\begin{theorem}\label{theorem:hausdorff_p_subadditive}
    Fix $1\leq p<\infty$. If $A,B\subset\mathbb{R}^n$ are nonempty compact sets, then 
    \begin{align*}
        d_0(A\oplus_pB)^p\leq d_0(A)^p+d_0(B)^p.
    \end{align*}
    In the case of $p=\infty$, we have
    \begin{align*}
        d_0(A\oplus_\infty B)\leq \max\{d_0(A),d_0(B)\}.
    \end{align*}
\end{theorem}

\begin{proof}
    By definition of $d_0$, we have
$\textup{conv}(0,A)\subset A+d_0(A)B_2^n$ and 
$\textup{conv}(0,B)\subset B+d_0(B) B_2^n$.
 Then
    \begin{align*}
        \textup{conv}(0,A\oplus _p B)&\subset \textup{conv}(0,A)\oplus _p \textup{conv}(0,B)\\
        &\subset (A+d_0(A)B_2^n)\oplus_p (B+d_0(B)B_2^n)\\
        &\subset A\oplus_p B+\|(d_0(A),d_0(B))\|_pB_2^n,
    \end{align*}
    and the desired bound follows.
\end{proof}
We now show that the power $p$ in Theorem \ref{theorem:hausdorff_p_subadditive} is optimal if $n\ge3$, adapting the counterexample \cite[Theorem 7.3]{FMMZ-18} to the setting of the $L_p$-sum.
\begin{prop}\label{proposition:sharp_dimension_larger_than_3_hausdorff_lp_sum}
    Let $n\geq 3$ be an integer, fix $1< p<\infty$ and choose $r>p$. Then for any $\alpha>0$, there exist compact sets $A_{\alpha},B_{\alpha}\subset\mathbb{R}^n$ such that $d_0(A_{\alpha})=d_0(B_{\alpha})=\alpha$,
    and 
    \begin{align*}
        d_0(A_{\alpha}\oplus_p B_{\alpha})>(d_0(A_{\alpha})^r+d_0(B_{\alpha})^r)^{1/r}.
    \end{align*}
    It follows that $d_0^r$ is subadditive with respect to $L_p$ summation of compact sets in $\mathbb{R}^n$ if and only if $r\leq p$. Moreover, if $r=p$, and if $\varepsilon>0$, then there exist compact sets $A_{\varepsilon,\alpha},B_{\varepsilon,\alpha}\subset\mathbb{R}^n$ such that $d_0(A_{\varepsilon,\alpha})=d_0(B_{\varepsilon,\alpha})=\alpha$, and 
    \begin{align*}
        d_0(A_{\varepsilon,\alpha}\oplus_p B_{\varepsilon,\alpha})\geq (1-\varepsilon)(d_0(A_{\varepsilon,\alpha})^p+d_0(B_{\varepsilon,\alpha})^p)^{1/p},
    \end{align*}
    showing that the inequality derived in Theorem \ref{theorem:hausdorff_p_subadditive} is sharp.
\end{prop}
\begin{proof}
We draw this specific example from the paper \cite{FMMZ-18}.

Let $\theta>0$ be a parameter to be specified later. For $\epsilon\in \{\pm 1\}$, denote by $A_{\epsilon}$ the segment $A_{\epsilon}=[0,(1, 0, \epsilon \theta)^\top]$, and by $B_{\epsilon}$ the segment $B_{\epsilon}=[0,(1, \epsilon \theta, 0)^\top]$. Let $A=A_{-1}\cup A_1$ and $B=B_1\cup B_{-1}$. In other words, if $\rho\in O(3)$ denotes the rotation of angle $\pi/2$ about the axis $\mathbb{R}e_1$, then $B=\rho(A)$. 
Hence, $\operatorname{conv}(A,B)$ is a pyramid with its apex at $0$ and a square base. Moreover, setting $F_{\varepsilon,\delta}:=A_{\varepsilon}\oplus_p B_{\delta}$  (when $\epsilon, \delta \in \{\pm 1\}$),  we may write for $p>1$ :

\[
    A\oplus_p B = \bigcup_{\epsilon, \delta\in \{-1,1\}} A_{\epsilon} \oplus_p B_{\delta} =: \bigcup_{\epsilon, \delta\in \{-1,1\}} F_{\epsilon, \delta} \, .
\]
Since $0\in A$, and by the symmetry of $A$ (with respect to $\R e_1$), one can easily check that
\[
    d_0(A) = d(A) = \sup_{z\in \operatorname{conv}(A)} d(z,A) = d\left(\begin{pmatrix} 1 \\ 0 \\ 0 \end{pmatrix}, A\right) = \frac{\theta}{\sqrt{1+\theta^2}},
\]
and since $B=\rho(A)$, we also have $d_0(B)=d(B)=\frac{\theta}{\sqrt{1+\theta^2}}$. 

Now let us compute $d_0(A\oplus_p B)$. Denote $a_1=(1,0,\theta)^\top$ and $b_1=(1,\theta,0)^\top$. Then, with $t=1/2$, we find that
\[
    z_+ := t^{1/q} a_1+(1-t)^{1/q} b_1 = \left(2^{1/p}, 2^{-1/q} \theta, 2^{-1/q} \theta\right)^\top \in F_{1,1} \subset A\oplus_p B,
\]
and similarly, 
$z_- := \left(2^{1/p}, -2^{-1/q} \theta, -2^{-1/q} \theta\right)^\top \in F_{-1,-1} \subset A\oplus_p B$. Thus, 
\[
    z = \frac{z_+ + z_-}{2} = (2^{1/p},0,0)^\top \in \operatorname{conv}(A\oplus_p B).
\]
Since $z=\rho(z)$ and since, for any $(i,j)\in \{\pm 1\}^2$, $F_{i,j}=\rho^k(F_{1,1})$ for some $k=k_{i,j}\in \{0,1,2,3\}$, we see that $d(z,F_{i,j})=d(z,F_{1,1})$ for all $(i,j)$. Thus,
\[
    d(z,A\oplus_p B) = \min_{i,j} d(z,F_{i,j}) = d(z,F_{1,1}).
\]
Since $F_{1,1}\subset \textup{span}(a_1,b_1)$, we see that $\operatorname{aff}(F_{1,1})=u^{\perp}$ with $u=\frac{1}{\sqrt{2+\theta^2}} (-\theta, 1,1)^\top$. Thus,
\begin{align*}
    P_{\R u}(z) &= \langle z,u\rangle u = -\frac{2^{1/p}\theta}{2+\theta^2} \begin{pmatrix} -\theta \\ 1 \\ 1 \end{pmatrix}, \\
    P_{u^{\perp}}(z) &= z - P_{\R u}(z) = \begin{pmatrix} 2^{1/p} \\ 0 \\ 0 \end{pmatrix} + \frac{2^{1/p}\theta}{2+\theta^2} \begin{pmatrix} -\theta \\ 1 \\ 1 \end{pmatrix} = \frac{2^{1/p}}{2+\theta^2} \begin{pmatrix} 2 \\ \theta \\ \theta \end{pmatrix}.
\end{align*}
Denote $p_{1,1}=P_{u^{\perp}}(z)=\frac{2^{1/p}}{2+\theta^2}(2,\theta,\theta)^\top$. We claim that $p_{1,1}\in F_{1,1}$. Indeed, set $\lambda=\frac{2}{2+\theta^2} \in [0,1]$. Since $F_{1,1}=[0,a_1]\oplus_p [0,b_1]$, we can verify this by checking:
\[
    p_{1,1} = \frac{2^{1/p}}{2+\theta^2} \begin{pmatrix} 2 \\ \theta \\ \theta \end{pmatrix} = 2^{-1/q}\frac{2}{2+\theta^2} \begin{pmatrix} 2 \\ \theta \\ \theta \end{pmatrix} = 2^{-1/q} \lambda \begin{pmatrix} 2 \\ \theta \\ \theta \end{pmatrix} = 2^{-1/q}(\lambda a_1+\lambda b_1) \in F_{1,1}.
\]
As a result, $d(z,F_{1,1})=d(z,u^{\perp})=|z-p_{1,1}|=|\langle z,u\rangle|=\frac{2^{1/p} \theta}{\sqrt{2+\theta^2}}$. Hence,
\[
    d(z,A\oplus_p B) = d\left(z, \bigcup_{i,j} F_{i,j}\right) = d(z,F_{1,1}) = \frac{2^{1/p} \theta}{\sqrt{2+\theta^2}}.
\]
If $r>p$, then we have
\begin{align*}
    \frac{d(A\oplus_p B)}{(d(A)^r+d(B)^r)^{1/r}}\geq 2^{1/p-1/r}\sqrt{\frac{1+\theta^2}{2+\theta^2}}.
\end{align*}
The above square root converges to $1$ as $\theta\rightarrow\infty$. Because $r>p$, we have $2^{1/r-1/p}<1$, so there exists $\theta$ so large that the square root term is strictly larger than $2^{1/r-1/p}$, and therefore the ratio is larger than $1$, giving us the desired strict inequality. If $r=p$, then $2^{1/p-1/r}=1$. In this case, just choose $\theta$ so large that the square root term is larger than $(1-\varepsilon)$, and we have what we want. To make sure $d(A)=d(B)=\alpha$, dilate the sets by an appropriate constant.
\end{proof}

\begin{theorem}\label{theorem:strong_fractional_subadditivity}
    Let $1< p\leq \infty$ and let $A_1,\dots,A_k$ be nonempty compact sets in $\mathbb{R}^n$ with $k\geq n+1$. Then
    \begin{align*}
        d_0\left(\bigoplus_{i\in [k]}^{(p)}A_i\right)\leq \max_{I\subset [k]:|I|\leq n}\min_{i\in [k]\backslash I}d_0\left(\bigoplus_{[k]\backslash \{i\}}^{(p)}A_j\right).
    \end{align*}
\end{theorem}

\begin{proof}
    Let $x\in \textup{conv}(0,\bigoplus_{i\in [k]}^{(p)}A_i)$. Then $x\in \bigoplus_{i\in [k]}^{(p)}\textup{conv}(0,A_i)$ and we can apply Lemma~\ref{shapley_folkman_LP} to find a set $I\subset [k]$ containing at most $n$ elements such that
    \begin{align*}
        x\in \left(\bigoplus_{i\in I}^{(p)}\textup{conv}(0,A_i)\right)\oplus_p
        \left(\bigoplus_{i\in [k]\backslash I}^{(p)}(\{0\}\cup A_i)\right).
    \end{align*}
    Choose $i_0\in [k]\backslash I$. We have
        $x\in \bigoplus_{i\in [k]\backslash \{i_0\}}^{(p)}\textup{conv}(0,A_i)\oplus_p (\{0\}\cup A_{i_0})$.
    Then we write
        $x=\sum_{i\in [k]\backslash \{i_0\}}\lambda_i^{1/q}x_i+\lambda_{i_0}^{1/q}x_{i_0}$,
    where $\lambda_i\geq 0$, $\sum_{i\in [k]}\lambda_i=1$, $x_i\in \textup{conv}(0,A_i)$ for $i\in [k]\backslash \{i_0\}$, and $x_{i_0}\in \{0\}\cup A_{i_0}$. Moreover, any $y\in \bigoplus_{i\in [k]}^{(p)}A_i$ is of the form
    \begin{align*}
        y=\sum_{i\in [k]\backslash \{i_0\}}\tilde{\lambda}_i^{1/q}y_i+\tilde{\lambda}_{i_0}^{1/q}y_{i_0},
    \end{align*}
    where $\tilde{\lambda}_i\geq 0$, $\sum_{i\in [k]}\tilde{\lambda}_i=1$, and $y_i\in A_i$ for all $i$. Suppose that $x_{i_0}\in A_{i_0}$. If $\lambda_{i_0}=1$, then $x=x_{i_0}\in A_{i_0}$, and $d(x,\bigoplus_{i\in [k]}^{(p)}A_i)=0$. If $\lambda_{i_0}<1$, then by restricting to the elements $y$ for which $\tilde{\lambda}_{i_0}=\lambda_{i_0}$ and $y_{i_0}=x_{i_0}$, we have
    \begin{align*}
        &d\left(x,\bigoplus_{i\in [k]}^{(p)}A_i\right)=\inf_{y\in \bigoplus_{i\in [k]}^{(p)} A_i}\|x-y\|_2\\
        &=\inf_{\substack{\sum_{i\in [k]}\tilde{\lambda}_i=1\\ \tilde{\lambda}_i\geq 0}}\left\|\sum_{i\in [k]\backslash \{i_0\}}\lambda_i^{1/q}x_i-\sum_{i\in [k]\backslash \{i_0\}}\tilde{\lambda}_i^{1/q}y_i+\lambda_{i_0}^{1/q}x_{i_0}-\tilde{\lambda}_{i_0}^{1/q}y_{i_0}\right\|_2\\
        &\leq (1-\lambda_{i_0})^{1/q}\inf_{\substack{\sum_{i\in [k]\backslash \{i_0\}}\tilde{\lambda}_i=1-\lambda_{i_0}\\ \tilde{\lambda}_i\geq 0}}\left\|\sum_{i\in [k]\backslash \{i_0\}}\left(\frac{\lambda_i}{1-\lambda_{i_0}}\right)^{1/q}x_i-\sum_{i\in [k]\backslash \{i_0\}}\left(\frac{\tilde{\lambda}_i}{1-\lambda_{i_0}}\right)^{1/q}y_i\right\|_2.
    \end{align*}
    Now, set $z=\sum_{i\in [k]\backslash \{i_0\}}(\lambda_i(1-\lambda_{i_0})^{-1})^{1/q}x_i$, which is an element of $\bigoplus_{i\in [k]\backslash \{i_0\}}^{(p)}\textup{conv}(0,A_i)$, and note that the above infimum is in fact over all $y\in \bigoplus_{i\in [k]\backslash \{i_0\}}^{(p)}A_i$. The above bound therefore implies
    \begin{align}\label{eq:leave_one_out_bound}
        d\left(x,\bigoplus_{i\in [k]}^{(p)}A_i\right)\leq d\left(z,\bigoplus_{i\in [k]\backslash \{i_0\}}^{(p)}A_i\right)\leq d_0\left(\bigoplus_{i\in [k]\backslash \{i_0\}}A_i\right).
    \end{align}
    If $x_{i_0}=0$, then restrict to those $y$ for which $\tilde{\lambda}_{i_0}=0$ to get
    \begin{align*}
        d\left(x,\bigoplus_{i\in [k]}^{(p)}A_i\right)&=\inf_{\substack{\sum_{i\in [k]}\tilde{\lambda}_i=1\\ \tilde{\lambda}_i\geq 0}}\left\|\sum_{i\in [k]\backslash \{i_0\}}\lambda_i^{1/q}x_i-\sum_{i\in [k]\backslash \{i_0\}}\tilde{\lambda}_i^{1/q}y_i-\tilde{\lambda}_{i_0}^{1/q}y_{i_0}\right\|_2\\
        &\leq \inf_{\substack{\sum_{i\in [k]\backslash \{i_0\}}\tilde{\lambda}_i=1\\ \tilde{\lambda}_i\geq 0}}\left\|\sum_{i\in [k]\backslash \{i_0\}}\lambda_i^{1/q}x_i-\sum_{i\in [k]\backslash \{i_0\}}\tilde{\lambda}_i^{1/q}y_i\right\|_2\\
        &\leq d_0\left(\bigoplus_{i\in [k]\backslash \{i_0\}}^{(p)}A_i\right),
    \end{align*}
    where for the last inequality we used the same technique as above. Note that this is possible because the element $z=\sum_{[k]\backslash \{i_0\}}\lambda_i^{1/q}x_i$ belongs to $\textup{conv}(0,\bigoplus_{i\in [k]\backslash \{i_0\}}A_i)$. Then the bound \eqref{eq:leave_one_out_bound} holds in any case for any $i_0\in [k]\backslash I$, so that 
    \begin{align*}
        d\left(x,\bigoplus_{i\in [k]}^{(p)}A_i\right)\leq \min_{i\in [k]\backslash I}d_0\left(\bigoplus_{i\in [k]\backslash \{i\}}^{(p)}A_i\right)
    \end{align*}
    Take the maximum over all sets $I\subset [k]$ of cardinality of at most $n$ to get
    \begin{align*}
        d\left(x,\bigoplus_{i\in [k]}^{(p)}A_i\right)\leq \max_{I\subset [k]:|I|\leq n}\min_{i\in [k]\backslash I}d_0\left(\bigoplus_{i\in [k]\backslash \{i\}}^{(p)}A_i\right)
    \end{align*}
    The result follows by taking the supremum over all $x\in \textup{conv}(0,\bigoplus_{i\in [k]}^{(p)}A_i)$.
\end{proof}

\begin{remark}
    Just as noted in the survey \cite{FMMZ-18}, if we consider the special case where $A_i=A$ for all $i\in [k]$, then we have the bound
    \begin{align*}
        d_0\left(\frac{1}{k^{1/p}}A^{\oplus_pk}\right)\leq \left(\frac{k-1}{k}\right)^{1/p}d_0\left(\frac{1}{(k-1)^{1/p}}A^{\oplus_p(k-1)}\right)\leq \left(\frac{n}{k}\right)^{1/p}d_0\left(\frac{1}{n^{1/p}}A^{\oplus_pn}\right).
    \end{align*}
    This gives us another proof of Theorem \ref{thm:Firey_convergence}.
\end{remark}

\section{Effective standard deviation}
Recall the effective standard deviation, denoted by $v$, which was defined in Section \ref{section:effective_standard_deviation}. We will prove the following generalization of Cassels' bound \eqref{eq:square_subadditive_cassels}.
\Pcassel*

\begin{proof}
    Again if $x\in \text{conv}(A)$, we denote $v_A(x)^2:=\inf\{V(X) : \text{supp}(X)\subset A, \mathbb EX=x\}$.
    Since $0\in \text{conv}(A)\cap \text{conv}(B)$, we know by Lemma \ref{lemma:conv0p},
    that $\text{conv}(A\oplus_p B)=\text{conv}(A)\oplus_p \text{conv}(B).$ Thus if $z\in \text{conv}(A\oplus_p B)$, we can find $x_1\in \text{conv}(A)$, $x_2\in \text{conv}(B)$, and $t_1\in [0,1]$, such that $z=t_1^{1/q} x_1+t_2^{1/q} x_2$, where $t_2=1-t_1$.

    If $\sum \lambda_i a_i=x_1$ and $\sum \mu_j b_j=x_2$ are two arbitrary (finite) representations of, respectively, $x_1$ and $x_2$ (i.e. both sums are finite, $\lambda_i >0$, $\mu_j>0$, $a_i\in A$, $b_j\in B$, and $\sum \lambda_i=\sum \mu_j=1$), then a (finite) representation of $z\in \text{conv}(A\oplus_p B)$ is given by:
    \[z=t_1^{1/q} x_1+t_2^{1/q} x_2=\sum_{i,j} \lambda_i \mu_j (t_1^{1/q} a_i+t_2^{1/q} b_j) \, .\]
    Therefore,
    \begin{align*}
    v(z,A\oplus_p B)^2 &\leq \sum_{i,j} \lambda_i \mu_j \big | t_1^{1/q} a_i+t_2^{1/q} b_j-z\big |^2=\sum_{i,j} \lambda_i \mu_j \big |t_1^{1/q}(a_i-x_1)+t_2^{1/q}(b_j-x_2)\big |^2 \\
    &= t_1^{2/q}\sum \lambda_i |a_i-x_1|^2 + t_2^{2/q} \sum \mu_j |b_j-x_2|^2.
    \end{align*}
    If $q>2$, then taking the infimum on the right-hand side, i.e., the infimum over the possible representations $x_1=\sum \lambda_i a_i \in \text{conv}(A)$; $x_2=\sum \mu_j b_j \in \text{conv}(B)$, one gets
    \begin{align*}
        &v(z,A\oplus_p B)^2 \leq t_1^{2/q} v(x_1,A)^2+t_2^{2/q} v(x_2,B)^2 \\
        &\leq (t_1+t_2)^{q/2}\cdot \left( v(x_1,A)^{\frac{2q}{q-2}} +v(x_2,B)^{\frac{2q}{q-2}}\right)^{\frac{q-2}{q}}=\left( v(x_1,A)^{\frac{2q}{q-2}} +v(x_2,B)^{\frac{2q}{q-2}}\right)^{\frac{q-2}{q}} 
 \end{align*}
    where we used H\"older's inequality, and $t_1+t_2=1$. If $1\leq q\leq 2$, then by convexity, we have 
    \begin{align*}
        v(z,A\oplus_p B)^2\leq t_1^{2/q}v(x_1,A)^2+t_2^{2/q}v(x_2,B)^2\leq \max\{v(x_1,A)^2,v(x_2,B)^2\}.
    \end{align*}
    Then taking the supremum over $z\in \text{conv}(A\oplus_p B)$ yields the claimed subadditivity.

    Next, for the planar case $n=2$, we explain why the maximum bound is true. First, since $0\in \textup{conv}(A)\cap \textup{conv}(B)$, Lemma \ref{lemma:conv0p} applies here with $\textup{conv}(A)=\textup{conv}(0,A)$ and $\textup{conv}(B)=\textup{conv}(0,B)$. Since $d(S)\leq v(S)$ holds for any set $S$, Lemma \ref{lemma:parallelogram_lemma} is valid if we use $v$ in place of $d$. Since $v(S)=d(S)$ when $S$ is the set of endpoints of an interval, Lemma \ref{lemma:firey_parallelogram_estimate} is valid for the distance function $v(\cdot,S)$ in place of $d(\cdot,S)$. Thus, recalling the formula
    \begin{align*}
        v(S)=\sup_{x\in \textup{conv}(S)}v(x,S),
    \end{align*}
    if we follow the proof strategy of Theorem \ref{theorem:lp_hausdorff_dimension_2_max_bound}, then we obtain the desired maximum bound for $v$ with the same example of equality. 

    For the case $n\geq 3$, let us prove the optimality of the bounds. When $2\leq p\leq \infty$, we use the exact same equality example from the equality case of Theorem \ref{theorem:lp_hausdorff_dimension_2_max_bound}. For the case $1< p<2$, choose positive numbers $\alpha_0$ and $\beta_0$ that satisfy
    \begin{align*}
        \alpha_0^q=\frac{a^s}{a^s+b^s},\qquad \beta_0^q=\frac{b^s}{a^s+b^s},
    \end{align*}
    where $s=2p/(2-p)$. Then 
    \begin{align*}
        \sqrt{a^2\alpha_0^2+b^2\beta_0^2}=(a^s+b^s)^{1/s}.
    \end{align*}
    Choose $\lambda>0$ so small that 
    \begin{align*}
        h_A:=\lambda\alpha_0^{q-1}\leq a,\qquad h_B:=\lambda\beta_0^{q-1}\leq b,
    \end{align*}
    and set 
    \begin{align*}
        &A:=\{0,h_Ae_3+ae_1,h_Ae_3-ae_1\},\qquad B:=\{0,h_Be_3+be_2,h_Be_3-be_2\}.
    \end{align*}
    It is straightforward to check that $v(A)=a$ and $v(B)=b$. A general point of the sumset $A\oplus_p B$ is of the form $u x+vy$ with $x\in A$, $y\in B$, and $u^q+v^q=1$. By H\"older's inequality, the $e_3$ component of such a point is at most
    \begin{align*}
        uh_A+vh_B\leq (h_A^p+h_B^p)^{1/p}=\lambda(\alpha_0^q+\beta_0^q)^{1/p}=\lambda.
    \end{align*}
    Therefore, the $e_3$ coordinate is of the form $x_3e_3$, with $x_3\leq \lambda$. Moreover, equality in H\"older's inequality is attained uniquely for $(u,v)=(\alpha_0,\beta_0)$. It follows that $A\oplus_p B$ intersects the plane $x_3=\lambda$ at exactly the four points $\lambda e_3\pm \alpha_0 ae_1\pm \beta_0 b e_2$. These points are the vertices of a rectangle with center $c:=\lambda e_3$.

    Since every point of $A\oplus_p B$ has $e_3$ component with $x_3\leq \lambda$, any convex combination of points from $A\oplus_p B$ that gives $c$ must use only points with $x_3=\lambda$. We showed above that the only such points are the vertices of the rectangle. Write $c=\sum p_iv_i$, where the $p_i$ are positive, and sum to $1$, and the $v_i$ are vertices of the mentioned rectangle above. Since $c$ is the center of the rectangle, $c$ is at an equal distance, $\ell$, from each vertex. Then, 
    \begin{align*}
        \sum p_i\|v_i-c\|_2^2=\ell^2.
    \end{align*}
    Elementary geometry gives
    \begin{align*}
        \ell=\sqrt{a^2\alpha_0^2+b^2\beta_0^2}=(a^s+b^s)^{1/s}.
    \end{align*}
    Putting everything together gives
    \begin{align*}
        v(c,A\oplus_p B)=\ell=(a^s+b^s)^{1/s},
    \end{align*}
    which verifies the equality.
\end{proof}

\begin{remark}
Let $A,B$ be two compact sets, and consider $A\oplus_p B$
for $p>1$ close to $1$. As $p\to 1$, we have, in the
Hausdorff metric,
\[
A\oplus_p B \longrightarrow
\left(\bigcup_{a\in A}[a,a+B]\right)
\cup
\left(\bigcup_{b\in B}[b,b+A]\right).
\]
In particular, $A+B$ is contained in the limit set, but
the inclusion may be strict.

Note that, even when both $A$ and $B$ contain the origin,
$A\oplus_p B$ need not converge to $A+B$ as $p\to 1$.
Thus, although $\frac{2p}{2-p}=2$ when $p=1$, one cannot
directly recover Cassels' result from
Theorem~\ref{prop:pcassel} (by letting $p\to 1$). To illustrate this, take $n=2$ and $A=\{0,a\}$, $B=\{0,b\},$ with $a=e_1, b=e_2+\frac{1}{4} e_1.$
Then $A+B=\{0,a,b,a+b\},$
whereas the Hausdorff limit of $A\oplus_p B$ as $p\to 1$
is
\[
D=[0,a]\cup[0,b]\cup[a,a+b]\cup[b,a+b],
\]
which is the boundary of the parallelogram
$\operatorname{conv}(A+B)$.

Applying Theorem~\ref{prop:pcassel} with
$p=1+\varepsilon$, where $0<\varepsilon<1$, gives
\[
v^{2\alpha_\varepsilon}(A\oplus_{1+\varepsilon}B)
\leq
v^{2\alpha_\varepsilon}(A)
+
v^{2\alpha_\varepsilon}(B),
\qquad
\alpha_\varepsilon=\frac{1+\varepsilon}{1-\varepsilon}
\longrightarrow 1.
\]
Letting $\varepsilon\to 0$, we obtain $v^2(D)\leq v^2(A)+v^2(B).$
And, since
$A+B\subset D$ and $\operatorname{conv}(A+B)=\operatorname{conv}(D)$,
we also have $v^2(D)\leq v^2(A+B)$, and for our example, one can check that this inequality is strict : it can be computed that $v^2(A+B)=\frac{425}{1024}$ while $v^2(D)=\frac{4}{17}$.
Thus since the inequality obtained by passing to the limit
$p\to 1$ involves the set $D$, rather than $A+B$, it is strictly weaker than Cassels' subadditivity inequality \eqref{eq:square_subadditive_cassels}, in
this example. In other words, the $p=1$ case cannot
be recovered directly from Theorem~\ref{prop:pcassel}
by letting $p\to 1$.
\end{remark}
We now modify the proof of \cite[Theorem 6.2]{FMMZ-18} to get the following strong subadditivity bound for $L_p$-sums.
\begin{theorem}\label{proposition:strong_fractional_subadditivity_v}
    Fix $1<p\leq \infty$. For an integer $k\geq n+1$, let $A_1,\dots,A_k\subset\mathbb{R}^n$ be compact sets  such that $0\in \textup{conv}(A_i)$ for all $i$. Then
    \begin{align*}
        v\left(\bigoplus_{i\in [k]}^{(p)}A_i\right)\leq \max_{I\subset [k]:|I|\leq n}\min_{i\in [k]\backslash I}v\left(\bigoplus_{j\in [k]\backslash \{i\}}^{(p)}A_j\right).
    \end{align*}
\end{theorem}

\begin{proof}
    Let $x\in \textup{conv}(\bigoplus_{i\in [k]}^{(p)}A_i)$. By Lemma \ref{lemma:conv0p}, $x\in \bigoplus_{i\in [k]}^{(p)}\textup{conv}(A_i)$. By Lemma \ref{shapley_folkman_LP}, there exists a set $I\subset [k]$ of cardinality at most $n$ such that 
    \begin{align*}
        x\in \left(\bigoplus_{i\in I}^{(p)}\textup{conv}(A_i)\right)\oplus_p\left(\bigoplus_{i\in [k]\backslash I}^{(p)}A_i\right).
    \end{align*}
    Choose $i_0\in [k]\backslash I$. Using Lemma \ref{shapley_folkman_LP} again, write
    \begin{align*}
        x\in \textup{conv}\left(\bigoplus_{i\in [k]\backslash \{i_0\}}^{(p)}A_i\right)\oplus_p A_{i_0}.
    \end{align*}
    There exists $0\leq t \leq 1$ such that $x=t^{1/q}z+(1-t)^{1/q}a_0$, where $z\in \textup{conv}(\bigoplus_{i \in [k]\setminus \{i_0\}}^{(p)}A_i)$ and $a_0\in A_{i_0}$. Then
    \begin{align*}
       & v\left(z,\bigoplus_{i\in [k]\backslash \{i_0\}}^{(p)}A_i\right)^2=\inf\left\{\sum_{m}p_m|\alpha_m-z|^2:z=\sum_{m}p_m\alpha_m,\sum_{m}p_m=1,\alpha_m\in \bigoplus_{i\in [k]\backslash\{i_0\}}^{(p)}A_i\right\}\\
        &\geq \inf\left\{\sum_{m}p_m|t^{1/q}\alpha_m-t^{1/q}z|^2:z=\sum_{m}p_m\alpha_m,\sum_{m}p_m=1,\alpha_m\in \bigoplus_{i\in [k]\backslash\{i_0\}}^{(p)}A_i\right\}
    \end{align*}
    Set $\alpha_m^{*}:=t^{1/q}\alpha_m+(1-t)^{1/q}a_0$, and note that $x=\sum_{m}p_m\alpha_m^*$. It is also clear from the definition of the $L_p$ sum that $\alpha_m^*\in \bigoplus_{i\in [k]}^{(p)}A_i$. By adding and subtracting $(1-t)^{1/q}a_0$, the above infimum is at least
    \begin{align*}
        \inf\left\{\sum_{m}p_m|\alpha_m^*-x|^2:x=\sum_{m}p_m\alpha_m^*,\sum_{m}p_m=1,\alpha_m\in \bigoplus_{i\in [k]}^{(p)}A_i\right\}=v\left(x,\bigoplus_{i\in [k]}^{(p)}A_i\right)^2
    \end{align*}
    Putting the above bounds together,
    \begin{align*}
        v\left(x,\bigoplus_{i\in [k]}^{(p)}A_i\right)^2\leq v\left(z,\bigoplus_{i\in [k]\backslash \{i_0\}}^{(p)}A_i\right)^2.
    \end{align*}
    Taking the supremum over $z\in \textup{conv}(\bigoplus_{i\in [k]\backslash\{i_0\}}^{(p)}A_i)$ and then minimizing over $i_0\in [k]$, we get
    \begin{align*}
        v\left(x,\bigoplus_{i\in [k]}^{(p)}A_i\right)^2\leq \min_{i\in [k]\backslash I}v\left(\bigoplus_{j\in [k]\backslash \{i\}}^{(p)}A_j\right)^2
    \end{align*}
    Then, as in the proof of Theorem \ref{theorem:strong_fractional_subadditivity}, take the max over all $I\subset [k]$ with $|I|\leq n$, and then over all $x\in \textup{conv}(\bigoplus_{i\in [k]}^{(p)}A_i)$ to arrive at the desired bound.
\end{proof}

\section{Volume of \texorpdfstring{$L_p$}{Lp}-sums}

In this section, we present the study of the volume of the $L_p$-sum of compact sets. For every compact set $A\subset\mathbb{R}^n$ and every $k\in\mathbb{N}$, the definition of $L_p$-sum gives
    \begin{equation} \label{eq:second-inclusion}
        A \subset \frac{A^{\oplus_p k}}{k^{1/p}} \subset \frac{(\conv (0,A))^{\oplus_p k}}{k^{1/p}} = \conv (0,A).
    \end{equation}
    Observe that the second inclusion may be strict. Indeed, the normalized sum need not contain the origin. However, by taking a convex hull of the origin and the sets, we obtain
    \[
        \conv (0,A) = \conv \left( 0, \frac{A^{\oplus_p k}}{k^{1/p}} \right).
    \]
    The theorem below shows that, for each fixed $1<p<\infty$, its volume converges to $|\operatorname{conv}(0,A)|$ as $k\to\infty$.
\ConvergenceVolume*

\begin{proof}
    We claim that for $k \geq 2n+2$, one has
    \begin{equation} \label{eq:Artem-inclusion}
        \left(1-\frac{n+1}{k}\right)^{1 / p} \operatorname{conv}(0, A) \setminus \left(\frac{n+1}{k}\right)^{1 / p} \operatorname{conv}(0, A) \subset \frac{1}{k^{1 / p}} A^{\oplus_p k} \subset \conv (0,A).
    \end{equation}
    The second inclusion follows from the observation above. We now prove the first inclusion. Let $x$ be a point in the leftmost set.  By Carath\'eodory's theorem, we write
    \[
        x =  \sum_{j=1}^{n+1} c_j a_j, \qquad a_j \in A, \quad c_j \geq 0, \quad \left(\frac{n+1}{k}\right)^{1 / p} < \sum_{j=1}^{n+1} c_j \leq \left(1-\frac{n+1}{k}\right)^{1 / p}.
    \]
    Set $u_j = k^{1/p}c_j$ and $S = \sum_{j=1}^{n+1} u_j$. By H\"older's inequality, we obtain
    \[
     \sum_{j=1}^{n+1} u_j ^q \geq \frac{S^q}{(n+1)^{q/p}}=\left(\frac{k}{n+1}\right)^{q/p}\left(\sum_{j=1}^{n+1}c_j\right)^q > 1.
    \]
    Choose
    $
    m_j=\left\lceil u_j S^{p-1}\right\rceil.
    $
    Then,
    \[
        \sum_{j=1}^{n+1} m_j \leq S^p+n+1 \leq k,
    \]
    where the last inequality follows from the upper bound on $\sum_{j=1}^{n+1} c_j$.  Increase some of the $m_j$'s so that their sum becomes exactly $k$, that is, $\sum_{j=1}^{n+1} m_j = k$.
    Moreover,
    $$
      \sum_{j: m_j>0} \frac{u_j^q}{m_j^{q/p} }\leq \sum_{j: u_j>0} \frac{u_j^q}{u_j^{q-1} S^{(p-1)(q-1)}}= \sum_{j=1}^{n+1} \frac{u_j}{S}=1.
    $$
    Applying Lemma \ref{lem:sum-of-ei} in dimension $n+1$, we get that $ u \in \{e_1,\ldots,e_{n+1}\}^{\oplus_p k}$. Define the linear map $T: \mathbb{R}^{n+1} \longrightarrow \mathbb{R}^n$ such that $T(e_j) = a_j$ for $1\leq j \leq n+1$. Hence,
\[
    k^{1 / p} x=T(u) \in T \left( \{e_1,\ldots,e_{n+1}\}^{\oplus_p k}\right) \subset \{a_1,\ldots,a_{n+1}\}^{\oplus_p k} \subset A ^{\oplus_p k},
\]
proving the first inclusion.
Taking volumes, we conclude that
\[
    0
    \leq
    \left|\operatorname{conv}(0,A)
    \setminus \frac{1}{k^{1 / p}} A^{\oplus_p k}\right|\leq
    \left[
        1-\left(1-\frac{n+1}{k}\right)^{n/p}
        +\left(\frac{n+1}{k}\right)^{n/p}
    \right]
    |\operatorname{conv}(0,A)|
    \longrightarrow 0.
\]
Together with
$\frac{1}{k^{1 / p}} A^{\oplus_p k}\subset\operatorname{conv}(0,A)$,
this yields $\lim_{k\to\infty}
    \left|\frac{1}{k^{1 / p}} A^{\oplus_p k}\right|
    =|\operatorname{conv}(0,A)|.$
\end{proof}

    \begin{remark} \label{rem:Dis-Hau}
        Hausdorff convergence alone does not imply volume convergence for arbitrary compact sets. Thus, the theorem above does not follow from the Hausdorff convergence established in Theorem~\ref{thm: hausd-convergence-2}. For example, the sets 
        $$
        E_k:=\left\{0, \frac{1}{k}, \ldots, \frac{k-1}{k}, 1\right\} \subset \mathbb{R}
        $$
        satisfy $ d_H\left(E_k,[0,1]\right)=\frac{1}{2 k} $, which converges to $ 0$ as $k \to \infty$, but $\left|E_k\right|=0.$
    \end{remark}
    \begin{remark}
    We note that the proof above in fact gives the rate of convergence of the volume for $1<p<\infty$, 
    \[
        |\operatorname{conv}(0, A)|-\left|\frac{1}{k^{1 / p}} A^{\oplus_p k}\right|=O\left(k^{-\min \{1, n / p\}}\right), \quad \text{as }k \to \infty.
    \]
    Moreover, the rate of convergence can be improved when $0 \in \conv (A)$. Indeed, Corollary \ref{coro:starshapedloneA} showed that $A^{\oplus_p k}$ is starshaped for $k \geq n+1$.
    By taking the radial hull on the inclusion \eqref{eq:Artem-inclusion}, we get that
    \begin{equation} \label{eq:inclusion-origin}
        \left(1-\frac{n+1}{k}\right)^{1 / p} \operatorname{conv}(0,A) \subset \frac{1}{k^{1 / p}} A^{\oplus_p k} \subset \operatorname{conv}(0,A), \qquad k \geq 2n+2.
    \end{equation}
    Therefore, the rate of convergence of the volume in this case is
    \[
        |\operatorname{conv}(0, A)|-\left|\frac{1}{k^{1 / p}} A^{\oplus_p k}\right|=O\left(k^{-1}\right).
    \]
\end{remark}

\subsection{The Bobkov, Madiman, Wang conjecture for \texorpdfstring{$L_p$}{Lp}-sums}\label{subsection:BMW_conjecture}

In the paper \cite{BMW-11}, Bobkov, Madiman, and Wang (BMW) studied a fractional generalization of the Brunn--Minkowski inequality. Their generalization is based on the concept of a fractional partition, which is defined as follows. Let $m\geq 2$ be an integer, let $\mathcal{G}$ be a nonempty collection of subsets of $[m]$, and let $\beta:\mathcal{G}\rightarrow [0,1]$ be a function of weights. The pair $(\mathcal{G},\beta)$ is called a fractional partition of $[m]$ if the weights satisfy the condition $\sum_{S\in \mathcal{G}:i\in S}\beta(S)=1$ for each $i\in [m]$.

Bobkov, Madiman, and Wang conjectured that if $A_1,\dots,A_m\subset\mathbb{R}^n$ are nonempty sets, and if $(\mathcal{G},\beta)$ is a fractional partition of $[m]$, then
\begin{align}\label{eq:fractional_brunn_minkowski}
    \left|\sum_{i=1}^{m}A_i\right|^{1/n}\geq \sum_{S\in \mathcal{G}}\beta(S)\left|\sum_{i\in S}A_i\right|^{1/n}.
\end{align}
If we set $m=2$, $\mathcal{G}=\{\{1\},\{2\}\}$, and $\beta(\{i\})=1$ for each $i$, then \eqref{eq:fractional_brunn_minkowski} is the Brunn--Minkowski inequality. The BMW conjecture is true if the sets $A_i$ are convex. In the case where the sets $A_i$ are not convex, the BMW conjecture was shown to be false \cite{FLZ-22} for $n\geq 7$, and true with the power $1/n$ removed \cite{BM-2024} with equality cases of the latter bound given in \cite{M2-2025}. The BMW conjecture is closely related to our study of convexification. We see that if $A_i=A$, $\mathcal{G}:=\{[m]\backslash \{i\}:i\in [m]\}$, and $\beta(S)=1/(m-1)$ for each $S\in \mathcal{G}$, then \eqref{eq:fractional_brunn_minkowski} reduces to the monotonicity of the volume of the $m$th set average $m^{-1}\sum_{i=1}^{m}A$. 

Our next theorem verifies an $L_p$-sum version of the BMW conjecture for the case where the sets $A_i$ are convex compact containing $0$. 
\begin{theorem}\label{fractional_superadditive_Lp_theorem}
    Let $A_1,\dots,A_m\subset \mathbb{R}^n$ be convex compact sets containing $0$, and let $(\mathcal{G},\beta)$ be a fractional partition of $[m]$. If $1<p<\infty$, then
    \begin{equation}\label{fractional_superadditive_Lp}
        \left|\bigoplus_{i\in [m]}^{(p)}A_i\right|^{p/n}\geq \sum_{S\in \mathcal{G}}\beta(S)\left|\bigoplus_{i\in S}^{(p)}A_i\right|^{p/n}.
    \end{equation}
\end{theorem}

\begin{proof}
    It is easy to see that
\begin{equation}\label{fractional_modify_sum}
        \bigoplus_{i\in[m]}^{(p)}A_i=\bigoplus_{S\in \mathcal{G}}^{(p)}\beta(S)^{1/p}\bigoplus_{i\in S}^{(p)}A_i.
    \end{equation}
    Now, use \eqref{fractional_modify_sum} and the $L_ p$-Brunn--Minkowski inequality \cite{LYZ-12} to obtain
    \begin{equation*}
        \begin{split}
            \left|\bigoplus_{i\in[m]}^{(p)}A_i\right|^{p/n}\geq\sum_{S\in\mathcal{G}}\left|\beta(S)^{1/p}\bigoplus_{i\in S}^{(p)}A_i \right|^{p/n}= \sum_{S\in\mathcal{G}}\beta(S)\left|\bigoplus_{i\in S}^{(p)}A_i\right|^{p/n}.
        \end{split}
    \end{equation*}
\end{proof}
Next, we present a counterexample motivated by \cite[Proposition 2]{FLZ-22}, which shows that convexity is necessary in Theorem \ref{fractional_superadditive_Lp_theorem}.
\begin{theorem}
    Let $1 \leq p<\infty$ and let $n \geq \lceil 7 p \rceil$. There exist starshaped, compact sets $A_1, A_2, A_3 \subset \mathbb{R}^n$ such that
$$
\left|A_1 \oplus_p A_2 \oplus_p A_3\right|^{p / n}<\frac{1}{2}\left(\left|A_1 \oplus_p A_2\right|^{p / n}+\left|A_1 \oplus_p A_3\right|^{p / n}+\left|A_2 \oplus_p A_3\right|^{p / n}\right) .
$$
\end{theorem}

\begin{proof}
    Let $m = \lfloor \frac{n}{2} \rfloor$, $Q_1 = [0,1]^{m} \times \{0\}^{n-m}$ and $Q_2 = \{0\}^m \times [0,1]^{n-m}$. Consider the sets
    \[
        A_1 = Q_1, \qquad A_2 = Q_2, \qquad A_3 = a^{1/p}Q_1 \cup b^{1/p} Q_2,
    \]
    where we select $a,b >0$ later. Then, $|A_1 \oplus_p A_2| = |Q_1 \oplus_p Q_2|$ and
    \[
    |A_1 \oplus_p A_3| = | (Q_1 \oplus_p a^{1/p}Q_1) \cup (Q_1 \oplus_p b^{1/p} Q_2) |.
    \]
    Since $Q_1 \oplus_p a^{1/p}Q_1 \subset \mathbb{R}^m \times\{0\}^{n-m}$, its volume are zero. Thus,
    \[
        |A_1 \oplus_p A_3| = | Q_1 \oplus_p b^{1/p} Q_2| =b^{(n-m)/p} |Q_1 \oplus_p Q_2|,
    \]
    where the last equality follows by applying twice the fact that for convex bodies $C$ and $D$ containing the origin and lying in complementary orthogonal subspaces of dimensions $m$ and $n-m$, respectively,
    \begin{equation} \label{eq:lp-on-complementary orthogonal subspaces}
    |C\oplus_p D|
    =
    \frac{
    \Gamma\left(1+\frac{m}{q}\right)
    \Gamma\left(1+\frac{n-m}{q}\right)
    }{
    \Gamma\left(1+\frac{n}{q}\right)
    }
    |C|_m\,|D|_{n-m}.
    \end{equation}
    Similarly, we have $|A_2 \oplus_p A_3| = a^{m/p} |Q_1 \oplus_p Q_2|$ and 
    \begin{align}
    |A_1 \oplus_p A_2 \oplus_p A_3| &= | (Q_1 \oplus_p Q_2 \oplus_p a^{1/p}Q_1) \cup (Q_1 \oplus_p Q_2 \oplus_p b^{1/p} Q_2)|
    \\
    & = \left| \left((1+a)^{1/p}Q_1 \oplus_p Q_2\right) \cup \left(Q_1 \oplus_p (1+b)^{1/p}Q_2\right) \right|
    \\
    &\leq \left| \left((1+a)^{1/p}Q_1 \oplus_p Q_2\right) \right| + \left| \left(Q_1 \oplus_p (1+b)^{1/p}Q_2\right) \right| - \left| Q_1 \oplus_p Q_2 \right|
    \\
    &= \left( (1+a)^{m/p} + (1+b)^{(n-m)/p} - 1 \right) |Q_1 \oplus_p Q_2|.
    \end{align}
    We recall the well-known Karamata's inequality: for any convex function $\varphi$, if $ (x,y)$ is majorized by $(u,v)$, that is, $u \leq x, y \leq v$ and $u+v = x+y$, then
    \[
        \varphi(x) + \varphi(y) \leq \varphi (u) +\varphi (v).
    \]
    Applying this with $ \varphi (t) = t^{n/(7p)}$, which is convex since $ n \geq 7p$, and 
    \[
        x = (1+a)^{7m/n}, \quad y = (1+b)^{7(n-m)/n}, \quad u =1, \quad v= x+y-1,
    \]
    we obtain
    \[
        \left( (1+a)^{m/p} + (1+b)^{(n-m)/p} - 1 \right)^{p/n} \leq \left((1+a)^{7m/n} +(1+b)^{7(n-m)/n} -1 \right)^{1/7}.
    \]
    Choose $ a ,b$ such that $a^{m/n}=b^{(n-m)/n} =2$. We need to verify that
    \begin{equation} \label{eq:desired-computation}
        \left(\left(1+2^{n/m}\right)^{7m/n} +\left(1+2^{n/(n-m)}\right)^{7(n-m)/n} -1 \right)^{1/7} < \frac{5}{2}.
    \end{equation}
    Note that $t \mapsto \left(1+2^{1 / t}\right)^{7 t}$ is  log-convex for $t>0$, and hence it is convex.
    Since $ m = \lfloor n/2 \rfloor$ and $ n \geq 7$, we get that $\frac{3}{7} \leq \frac{m}{n} \leq \frac{n-m}{n} \leq \frac{4}{7}$. Applying Karamata's inequality yields
    \[
        \left(1+2^{n / m}\right)^{7 m / n}+\left(1+2^{n /(n-m)}\right)^{7(n-m) / n} \leq\left(1+2^{7 / 3}\right)^3+\left(1+2^{7 / 4}\right)^4 < \left(\frac{5}{2}\right)^7+1,
    \]
    where the last inequality is a direct computation. The proof is complete since it implies \eqref{eq:desired-computation}.
\end{proof}

\subsection{Eventual monotonicity of volumes}

Recall from Section \ref{subsection:BMW_conjecture} that the BMW conjecture reduces to a question of monotonicity of the volume of Minkowski averages. Stated precisely, for a nonempty compact set $A \subset \mathbb{R}^n$, we ask whether the sequence
$$
\left\{\left|\frac{1}{k} \sum_{i=1}^{k}A\right|\right\}_{k \geq 1}
$$
 is nondecreasing in $k$. A negative answer for $n \geq 12$ was given in \cite[Theorem 3.4]{FMMZ-18}. However,  Fradelizi, L\'angi, and Zvavitch \cite{FLZ-22} showed that if $A$ is starshaped, then this sequence is nondecreasing for dimension $n=2,3$, and eventually nondecreasing for dimension $n\geq4$. We will establish an analogous result for the $L_p$-sum averages, showing that the sequence of volumes
$$
\left\{\left|\frac{1}{k^{1/p}} A^{\oplus_p k}\right|\right\}_{k \geq 1}
$$
is eventually nondecreasing.

We begin with a counterexample based on the same construction as in \cite{FMMZ-18}, recovering their result when $p=1$.
\begin{theorem} \label{thm:counterexample-12}
Let $1 \leq p<\infty$, and let $n \geq 
2\left\lceil\frac{p \ln 2}{\ln (9 / 8)}\right\rceil
$. There exists a starshaped set $A$ such that 
\[
    \left|\frac{1}{3^{1/p}}A^{\oplus_p 3}\right| < \left|\frac{1}{2^{1/p}}A^{\oplus_p 2}\right|.
\]
\end{theorem}

\begin{proof}
    Let $m = \lfloor \frac{n}{2} \rfloor$ and $A = Q_1 \cup Q_2$ where $Q_1 = [0,1]^{m} \times \{0\}^{n-m}$, $Q_2 = \{0\}^m \times [0,1]^{n-m}$.
    Observe that
    \[
        \frac{1}{2^{1/p}} A^{\oplus_p 2} = \frac{1}{2^{1/p}} \bigg((Q_1 \oplus_p Q_1) \cup (Q_2 \oplus_p Q_2) \cup (Q_1 \oplus_p Q_2) \bigg).
    \]
    Since $Q_1 \oplus_p Q_1 \subset \R^m \times \{0\}^{n-m}$ and $Q_2 \oplus_p Q_2 \subset \{ 0 \}^m \times \R^{n-m}$, their volumes are zero.
    Thus, $|\frac{1}{2^{1/p}}(A\oplus_p A)|_n = \frac{1}{2^{n/p}} |Q_1 \oplus_p Q_2|_n$. 
    Next, we compute the volume of the left-hand side,
    \begin{align}
        \left|\frac{1}{3^{1/p}}A^{\oplus_p 3}\right| 
        &= \frac{1}{3^{n/p}} \left| \left( Q_1\oplus_p Q_1 \oplus_p Q_2 \right)
        \cup \left(Q_1 \oplus_p Q_2 \oplus_p Q_2 \right) \right|
        \\
        &= \frac{1}{3^{n/p}} \left| (2^{1/p}Q_1 \oplus_p Q_2) \cup (Q_1 \oplus_p 2^{1/p} Q_2)  \right| 
        \\
        &\leq \frac{1}{3^{n/p}} \bigg( |2^{1/p}Q_1 \oplus_p Q_2| + |Q_1 \oplus_p 2^{1/p}Q_2| - |Q_1 \oplus_p Q_2| \bigg)
        \\
        &= \left(\frac{2^{m/p} + 2^{(n-m)/p} -1}{3^{n/p}} \right) |Q_1 \oplus_p Q_2|,
    \end{align}
    where the last equality used \eqref{eq:lp-on-complementary orthogonal subspaces}.
    It is enough to check that if $n \geq 2\left\lceil\frac{p \ln 2}{\ln (9 / 8)}\right\rceil$, then
    \[
        \frac{2^{m/p} + 2^{(n-m)/p} }{3^{n/p}} \leq \frac{1}{2^{n/p}}.
    \]
    Using that $ m = \lfloor \frac{n}{2} \rfloor$, we have that $\delta := n-2m \in \{0,1\}$. 
    In terms of $\delta$, we want to show that
    \begin{align}
        \left( \frac{2}{3} \right)^{\delta/p} +\left( \frac{4}{3} \right)^{\delta/p} \leq \left( \frac{9}{8} \right)^{m/p}.
    \end{align}
    Note that
    \[
        m = \left\lfloor \frac{n}{2} \right\rfloor \geq \left\lceil\frac{p \ln 2}{\ln (9 / 8)}\right\rceil \geq \frac{p \ln 2}{\ln (9 / 8)}.
    \]
    Then, $\left( \frac{9}{8} \right)^{m/p}\geq 2.$ On the other hand, we consider the function $\phi $ on $[0,1]$ defined by
    \[
        \phi (t) = \left( \frac{2}{3} \right)^{t} +\left( \frac{4}{3} \right)^{t}.
    \]
    Note that $\phi$ is convex with $\phi(0) = \phi(1) =2$, so we get that
    \[
        \left( \frac{2}{3} \right)^{\delta/p} +\left( \frac{4}{3} \right)^{\delta/p} \leq 2 \leq \left( \frac{9}{8} \right)^{m/p}. \qedhere
    \]
\end{proof}

Next, we establish the eventual monotonicity of the volumes of $L_p$-sum averages. We first recall a volume estimate from the proof of Lemma~2 in \cite{FLZ-22}.
\begin{lemma}[{\cite[Eq.~(10)]{FLZ-22}}]
\label{lem:aux-monotonic-vol-1bis}
Let $n\geq 2$ and let $B=\bigcup_{i=1}^n[0,e_i]$. If $k\geq\max\{2,(n-1)(n-2)\}$, and if $M\subset\R^n$ is a  compact set with  $B^{\oplus_1 k}\subset M\subset k\operatorname{conv}(B)$, then
\[
\left|\frac{M+B}{k+1}\right|
\geq
\theta_{n,k}\left|\frac{M}{k}\right|
+(1-\theta_{n,k})|\operatorname{conv}(B)|, \quad \text{ where} \quad 
\theta_{n,k}
=\frac{k^n}{(k-n+2)(k+1)^{n-1}}<1.
\]
\end{lemma}

\begin{lemma}
\label{lem:aux-monotonic-vol}
Let $n\geq 2$, $1<p<\infty$, and let
$B=\bigcup_{i=1}^n[0,e_i]$.
If $k \geq \max \{2,(n-1)(n-2)\}$, and $M \subset \mathbb{R}^n$ is compact with 
    $B^{\oplus_p k} \subset M \subset k^{1 / p} \operatorname{conv}(B)$, then
\begin{equation} \label{eq:first-aux-monotonic-vol}
\left|\frac{M\oplus_p B}{(k+1)^{1/p}}\right|
\geq
\left|\frac{M}{k^{1/p}}\right|.
\end{equation}
Equality holds if and only if
$M=k^{1/p}\operatorname{conv}(B)$.
Moreover, there exists a constant $C(n,k)>0$,
depending only on $n$ and $k$, such that
\begin{equation} \label{eq:second-aux-monotonic-vol}
|\operatorname{conv}(B)|
-\left|\frac{M}{k^{1/p}}\right|
\leq C(n,k)
\left(
\left|\frac{M\oplus_p B}{(k+1)^{1/p}}\right|
-\left|\frac{M}{k^{1/p}}\right|
\right).
\end{equation}
\end{lemma}
\begin{proof}
By definition of the $L_p$-sum,
\[
B^{\oplus_1 k}\subset k^{1/q}B^{\oplus_p k} \subset k^{1/q}M \subset k\operatorname{conv}(B),
\qquad
\frac{k^{1/q}M+B}{k+1}
\subset\frac{M\oplus_p B}{(k+1)^{1/p}}
\subset 
\operatorname{conv}(B).
\]
Applying Lemma \ref{lem:aux-monotonic-vol-1bis} to $k^{1/q}M$ yields
\[
\left|\frac{M\oplus_p B}{(k+1)^{1/p}}\right|
\geq \left|\frac{k^{1/q}M+B}{k+1}\right| \geq \theta_{n,k}\left|\frac{M}{k^{1/p}}\right|
 +(1-\theta_{n,k})|\operatorname{conv}(B)|.
\]
Consequently,
\[
\left|\frac{M\oplus_p B}{(k+1)^{1/p}}\right|
-\left|\frac{M}{k^{1/p}}\right|
\geq (1-\theta_{n,k})
\left(|\operatorname{conv}(B)|-\left|\frac{M}{k^{1/p}}\right|\right)
\geq 0.
\]
This proves both inequalities with $C(n,k)=(1-\theta_{n,k})^{-1}$.

Finally, suppose equality holds in \eqref{eq:first-aux-monotonic-vol}. Since $M$ is compact and
$M\subset k^{1/p}\operatorname{conv}(B)$, it follows from \eqref{eq:second-aux-monotonic-vol} that  $M=k^{1/p}\operatorname{conv}(B)$. Conversely, if $M=k^{1/p}\operatorname{conv}(B)$, then 
\[
    \left|\frac{M}{k^{1/p}}\right| = |\operatorname{conv}(B)| \geq \left| \frac{M \oplus_p B}{(k+1)^{1 / p}} \right|,
\]
where the last inequality follows from the inclusion above. Together with \eqref{eq:first-aux-monotonic-vol}, this gives equality.
\end{proof}

\begin{theorem} \label{thm:even-mono-vol}
    Let $ n \geq 2$, $ 1 < p < \infty$, and let $A $ be a compact starshaped set. Then, for every $k \geq \max \{2,(n-1)(n-2)\}$,
    \begin{equation} \label{eq: monotone-col}
        \left|\frac{A^{\oplus_p(k+1)}}{(k+1)^{1 / p}}\right| \geq\left|\frac{A^{\oplus_p k}}{k^{1 / p}}\right| .
    \end{equation}
    Equality holds if and only if $A$ lies in a proper linear subspace or $A^{\oplus_p k} =k^{1 / p} \operatorname{conv}(A) .$
\end{theorem}

\begin{proof}
    We may assume that $\Span (A) = \R^n$. Let $ \varepsilon >0$. Choose a finite set $A_0 \subset A$ such that $|\operatorname{conv}(A) \setminus \operatorname{conv}(0,A_0) | < \varepsilon$. Since $\frac{A^{\oplus_p k}}{k^{1 / p}} \subset \operatorname{conv}(A)$,
    we obtain
    \begin{equation} \label{eq: monotone-col-01}
        \left| \frac{A^{\oplus_p k}}{k^{1 / p}}  \right| - \varepsilon \leq \left| \frac{A^{\oplus_p k}}{k^{1 / p}} \cap \operatorname{conv}(0,A_0) \right|.
    \end{equation}
    Since $\operatorname{conv}(0,A_0)$ is a polytope, we can triangulate it as the union of pairwise disjoint full-dimensional simplices with one vertex at the origin. We denote these simplices by 
    $$
        T_j=\operatorname{conv}\left(0, v_{j, 1}, \ldots, v_{j, n}\right), \quad v_{j,i} \in A_0,\quad j=1, \ldots, m.
    $$
    For each $j$, set
    $$
    B_j=\bigcup_{i=1}^n\left[0, v_{j, i}\right], \quad M_j=A^{\oplus_p k} \cap k^{1 / p} T_j .
    $$
    Since $A$ is starshaped, $B_j \subset A$. Hence,
    $
    B_j^{\oplus_p k} \subset M_j \subset k^{1 / p} T_j .
    $
    After applying an invertible linear map, Lemma \ref{lem:aux-monotonic-vol} gives
    \[
        \left|T_j\right|-\left|\frac{M_j}{k^{1 / p}}\right| \leq C(n, k)\left(\left|\frac{M_j \oplus_p B_j}{(k+1)^{1 / p}}\right|-\left|\frac{M_j}{k^{1 / p}}\right|\right),
    \]
    where $C(n,k)$ is independent of $j$. Note that
    \[
        \frac{M_j \oplus_p B_j}{(k+1)^{1 / p}} \subset  \frac{A^{\oplus_p (k+1)}}{(k+1)^{1 / p}} \cap T_j.
    \]
    Summing over $j$, we obtain
    \begin{align}
    \left|\operatorname{conv}\left(0, A_0\right) \setminus \frac{A^{\oplus_p k}}{k^{1 / p}} \right|
    & =\sum_{j=1}^m\left(\left|T_j\right|-\left|\frac{M_j}{k^{1 / p}}\right|\right)   
    \\
    &\leq C(n, k)\left(\left|\frac{A^{\oplus_p(k+1)}}{(k+1)^{1 / p}}\right| -\left|\frac{A^{\oplus_p k}}{k^{1 / p}} \cap \operatorname{conv}\left(0, A_0\right)\right|\right) .
    \end{align}
    By \eqref{eq: monotone-col-01}, we get
    \begin{align}
     |\operatorname{conv}(A)|-\left|\frac{A^{\oplus_p k}}{k^{1 / p}}\right|
    &  \leq \varepsilon+\left|\operatorname{conv}\left(0, A_0\right) \setminus \frac{A^{\oplus_p k}}{k^{1 / p}}\right| \\
    &  \leq(1+C(n, k)) \varepsilon+C(n, k)\left(\left|\frac{A^{\oplus_p(k+1)}}{(k+1)^{1 / p}}\right|-\left|\frac{A^{\oplus_p k}}{k^{1 / p}}\right|\right) .
    \end{align}
    Letting $ \varepsilon \to 0$ gives
    \begin{align}
         |\operatorname{conv}(A)|-\left|\frac{A^{\oplus_p k}}{k^{1 / p}}\right|
        &  \leq C(n, k)\left(\left|\frac{A^{\oplus_p(k+1)}}{(k+1)^{1 / p}}\right|-\left|\frac{A^{\oplus_p k}}{k^{1 / p}}\right|\right) .
    \end{align}
    The left-hand side is nonnegative, proving the desired monotonicity.
    
    For the equality case, suppose that equality holds in \eqref{eq: monotone-col}. The preceding estimate gives
    \[
        \left|\operatorname{conv}\left(A\right) \setminus \frac{A^{\oplus_p k} }{k^{1 / p}}\right| = 0.
    \]
    Thus,
    \[
        \frac{A^{\oplus_p k} }{k^{1 / p}}=\operatorname{conv}(A).
    \]
    The reverse direction is immediate since 
    $$
        \frac{A^{\oplus_p(k+1)}}{(k+1)^{1 / p}} \subset \operatorname{conv}(A) = \frac{A^{\oplus_pk}}{(k)^{1 / p}} ,
    $$ 
    completing the proof.
\end{proof}
Next, we extend the theorem above to sets satisfying weaker assumptions. The two theorems have the same lower bound on $k$, except in dimension $ n =3$.
\MonotoneCol*

\begin{proof}
    We may assume that $\Span (A) = \R^n$. Recall that $A^* = \bigcup_{a \in A}[0, a]$ and, for $k\ge n$, define
    $$
        H_{k}=(A^*)^{\oplus_pk} \setminus A^{\oplus_pk}.
    $$
    Corollary \ref{cor:subtraction-starshaped} shows that $H_{k+1}$ is either empty or starshaped and is nonincreasing in $k$. Now, we apply Theorem \ref{thm:even-mono-vol},
    \begin{align}
    \left|\frac{A^{\oplus_p(k+1)}}{(k+1)^{1 / p}}\right|  &=\left|\frac{(A^*)^{\oplus_p(k+1)}}{(k+1)^{1 / p}}\right|-\frac{|H_{k+1}|}{(k+1)^{n / p}}  \geq\left|\frac{(A^*)^{\oplus_p k}}{k^{1 / p}}\right|-\frac{|H_{k+1}|}{(k+1)^{n / p}} 
    \\
    &\geq\left|\frac{(A^*)^{\oplus_p k}}{k^{1 / p}}\right|-\frac{|H_{k}|}{(k+1)^{n / p}}  
    \geq\left|\frac{(A^*)^{\oplus_p k}}{k^{1 / p}}\right|-\frac{|H_{k}|}{k^{n / p}}  
    =\left|\frac{A^{\oplus_p k}}{k^{1 / p}}\right| .
    \end{align}
    This proves the desired inequality.

    Suppose that equality holds in \eqref{eq: monotone-col-compact}. Then, all inequalities in the above chain are equalities. Equality in the first inequality, together with the equality case in Theorem \ref{thm:even-mono-vol}, gives
    \[
        (A^*)^{\oplus_p k}=k^{1 / p} \operatorname{conv}(A^*)=k^{1 / p} \operatorname{conv}(0, A) .
    \]
    Equality in the third inequality yields $|H_k| = 0$ and so $\left|A^{\oplus_p k}\right|=\left|(A^*)^{\oplus_p k}\right| $. Thus, $A^{\oplus_p k}=(A^*)^{\oplus_p k} = k^{1 / p} \operatorname{conv}(0, A)$. 
    
    The reverse direction is immediate since 
    $$
        \frac{A^{\oplus_p(k+1)}}{(k+1)^{1 / p}} \subset \operatorname{conv}(0,A) = \frac{A^{\oplus_pk}}{k^{1 / p}} ,
    $$ 
    completing the proof.
\end{proof}



\end{document}